\documentclass[11pt,twoside]{article}
\usepackage{amsmath, amsthm, amscd, amsfonts, amssymb, graphicx, color}

\date{}

\begin{document}

\centerline{}

\centerline {\Large{\bf Weaving continuous controlled $K$-$g$-fusion frames }}
\centerline {\Large{\bf  in Hilbert spaces }}

\newcommand{\mvec}[1]{\mbox{\bfseries\itshape #1}}
\centerline{}
\centerline{\textbf{Prasenjit Ghosh}}
\centerline{Department of Pure Mathematics, University of Calcutta,}
\centerline{35, Ballygunge Circular Road, Kolkata, 700019, West Bengal, India}
\centerline{e-mail: prasenjitpuremath@gmail.com}
\centerline{}
\centerline{\textbf{T. K. Samanta}}
\centerline{Department of Mathematics, Uluberia College,}
\centerline{Uluberia, Howrah, 711315,  West Bengal, India}
\centerline{e-mail: mumpu$_{-}$tapas5@yahoo.co.in}

\newtheorem{Theorem}{\quad Theorem}[section]

\newtheorem{definition}[Theorem]{\quad Definition}

\newtheorem{theorem}[Theorem]{\quad Theorem}

\newtheorem{remark}[Theorem]{\quad Remark}

\newtheorem{corollary}[Theorem]{\quad Corollary}

\newtheorem{note}[Theorem]{\quad Note}

\newtheorem{lemma}[Theorem]{\quad Lemma}

\newtheorem{example}[Theorem]{\quad Example}

\newtheorem{result}[Theorem]{\quad Result}
\newtheorem{conclusion}[Theorem]{\quad Conclusion}

\newtheorem{proposition}[Theorem]{\quad Proposition}

\begin{abstract}
\textbf{\emph{We introduce the notion of weaving continuous controlled $K$-$g$-fusion frame in Hilbert space.\,Some characterizations of weaving continuous controlled $K$-$g$-fusion frame have been presented.\,We extend some of the recent results of woven $K$-$g$-fusion frame and controlled $K$-$g$-fusion frame to woven continuous controlled $K$-$g$-fusion frame. Finally, a perturbation result of woven continuous controlled $K$-$g$-fusion frame has been studied.}}
\end{abstract}
\\
\\

\section{Introduction and Preliminaries}
 
\smallskip\hspace{.6 cm} 
Duffin and Schaeffer \cite{Duffin} introduced frame for Hilbert space to study some fundamental problems in non-harmonic Fourier series.\,Later on, after some decades, frame theory was popularized by Daubechies et al.\,\cite{Daubechies}.\,At present, frame theory has been widely used in signal and image processing, filter bank theory, coding and communications, system modeling and so on.

Let \,$H$\, be a separable Hilbert space associated with the inner product \,$\left<\,\cdot,\,\cdot\,\right>$.\,Frame for Hilbert space was defined as a sequence of basis-like elements in Hilbert space.\,A sequence \,$\left\{\,f_{\,i}\,\right\}_{i \,=\, 1}^{\infty} \,\subset\, H$\, is called a frame for \,$H$, if there exist positive constants \,$0 \,<\, A \,\leq\, B \,<\, \infty$\, such that
\[ A\; \|\,f\,\|^{\,2} \,\leq\, \sum\limits_{i \,=\, 1}^{\infty}\, \left|\ \left <\,f \,,\, f_{\,i} \, \right >\,\right|^{\,2} \,\leq\, B \,\|\,f\,\|^{\,2}\; \;\text{for all}\; \;f \,\in\, H.\]
The constants \,$A$\, and \,$B$\, are called lower and upper bounds, respectively.

Throughout this paper,\;$H$\; is considered to be a separable Hilbert space with associated inner product \,$\left <\,\cdot \,,\, \cdot\,\right>$\, and \,$\mathbb{H}$\, is the collection of all closed subspaces of \,$H$.\,$(\,X,\,\mu\,)$\, denotes abstract measure space with positive measure \,$\mu$.\,$I_{H}$\; is the identity operator on \,$H$.\,$\mathcal{B}\,(\,H_{\,1},\, H_{\,2}\,)$\; is a collection of all bounded linear operators from \,$H_{\,1} \,\text{to}\, H_{\,2}$.\,In particular \,$\mathcal{B}\,(\,H\,)$\, denotes the space of all bounded linear operators on \,$H$.\;For \,$S \,\in\, \mathcal{B}\,(\,H\,)$, we denote \,$\mathcal{N}\,(\,S\,)$\; and \,$\mathcal{R}\,(\,S\,)$\, for null space and range of \,$S$, respectively.\,Also, \,$P_{M} \,\in\, \mathcal{B}\,(\,H\,)$\; is the orthonormal projection onto a closed subspace \,$M \,\subset\, H$.\,The set \,$\mathcal{S}\,(\,H\,)$\; of all self-adjoint operators on \,$H$\; is a partially ordered set with respect to the partial order \,$\leq$\, which is defined as for \,$R,\,S \,\in\, \mathcal{S}\,(\,H\,)$ 
\[R \,\leq\, S \,\Leftrightarrow\, \left<\,R\,f,\, f\,\right> \,\leq\, \left<\,S\,f,\, f\,\right>\; \;\forall\; f \,\in\, H.\]

\;$\mathcal{G}\,\mathcal{B}\,(\,H\,)$\, denotes the set of all bounded linear operators which have bounded inverse.\,If \,$S,\, R \,\in\, \mathcal{G}\,\mathcal{B}\,(\,H\,)$, then \,$R^{\,\ast},\, R^{\,-\, 1}$\, and \,$S\,R$\, also belongs to \,$\mathcal{G}\,\mathcal{B}\,(\,H\,)$.\,A self-adjoint operator \,$U \,:\, H \,\to\, H$\, is called positive if \,$\left<\,U\,f \,,\,  f\,\right> \,\geq\, 0$\, for all \,$f \,\in\, H$.\;In notation, we can write \,$U \,\geq\, 0$.\;A self-adjoint operator \,$V \,:\, H \,\to\, H$\, is called a square root of \,$U$\, if \,$V^{\,2} \,=\, U$.\;If, in addition \,$V \,\geq\, 0$, then \,$V$\, is called positive square root of \,$U$\, and is denoted by \,$V \,=\, U^{1 \,/\, 2}$.\,The positive square root \,$V \,:\, H \,\to\, H$\, of an arbitrary positive self-adjoint operator \,$U \,:\, H \,\to\, H$\, exists and is unique.\;Further, the operator \,$V$\, commutes with every bounded linear operator on \,$H$\, which commutes with \,$U$.\;$\mathcal{G}\,\mathcal{B}^{\,+}\,(\,H\,)$\, is the set of all positive operators in \,$\mathcal{G}\,\mathcal{B}\,(\,H\,)$\, and \,$T,\, U$\, are invertible operators in \,$\mathcal{G}\,\mathcal{B}\,(\,H\,)$.\,For each \,$m \,>\, 1$, we define \,$[\,m\,] \,=\, \{\,1,\, 2,\, \cdots,\, m\,\}$. 

We present some theorems in operator theory which will be needed throughout this paper.

\begin{theorem}(\,Douglas' factorization theorem\,)\,{\cite{Douglas}}\label{th1}
Let \;$S,\, V \,\in\, \mathcal{B}\,(\,H\,)$.\,Then the following conditions are equivalent:
\begin{description}
\item[$(i)$]$\mathcal{R}\,(\,S\,) \,\subseteq\, \mathcal{R}\,(\,V\,)$.
\item[$(ii)$]\;\;$S\, S^{\,\ast} \,\leq\, \lambda^{\,2}\; V\,V^{\,\ast}$\; for some \,$\lambda \,>\, 0$.
\item[$(iii)$]$S \,=\, V\,W$\, for some bounded linear operator \,$W$\, on \,$H$.
\end{description}
\end{theorem}

\begin{theorem}\cite{Gavruta}\label{th1.01}
Let \,$M \,\subset\, H$\; be a closed subspace and \,$T \,\in\, \mathcal{B}\,(\,H\,)$.\;Then \,$P_{\,M}\, T^{\,\ast} \,=\, P_{\,M}\,T^{\,\ast}\, P_{\,\overline{T\,M}}$.\;If \,$T$\; is an unitary operator (\,i\,.\,e \,$T^{\,\ast}\, T \,=\, I_{H}$\,), then \,$P_{\,\overline{T\,M}}\;T \,=\, T\,P_{\,M}$.
\end{theorem}

\subsection{$K$-$g$-fusion frame}
\smallskip\hspace{.6 cm}
Construction of \,$K$-$g$-fusion frames and their dual were presented by Sadri and Rahimi \cite{Sadri} to generalize the theory of \,$K$-frame \cite{L}, fusion frame \cite{Kutyniok}, and \,$g$-frame \cite{Sun}.

\begin{definition}\cite{Sadri}
Let \,$\left\{\,W_{j}\,\right\}_{ j \,\in\, J}$\, be a collection of closed subspaces of \,$H$\; and \,$\left\{\,v_{j}\,\right\}_{ j \,\in\, J}$\; be a collection of positive weights, \,$\left\{\,H_{j}\,\right\}_{ j \,\in\, J}$\, be a sequence of Hilbert spaces.\,Suppose \,$\Lambda_{j} \,\in\, \mathcal{B}\,(\,H,\, H_{j}\,)$\; for each \,$j \,\in\, J$\, and \,$K \,\in\, \mathcal{B}\,(\,H\,)$.\;Then \,$\Lambda \,=\, \{\,\left(\,W_{j},\, \Lambda_{j},\, v_{j}\,\right)\,\}_{j \,\in\, J}$\; is called a $K$-$g$-fusion frame for \,$H$\; respect to \,$\left\{\,H_{j}\,\right\}_{j \,\in\, J}$\; if there exist constants \,$0 \,<\, A \,\leq\, B \,<\, \infty$\, such that
\[A \;\left \|\,K^{\,\ast}\,f \,\right \|^{\,2} \,\leq\, \sum\limits_{\,j \,\in\, J}\,v_{j}^{\,2}\, \left\|\,\Lambda_{j}\,P_{\,W_{j}}\,(\,f\,) \,\right\|^{\,2} \,\leq\, B \; \left\|\, f \, \right\|^{\,2}\; \;\forall\; f \,\in\, H.\]
The constants \,$A$\; and \,$B$\; are called the lower and upper bounds of $K$-$g$-fusion frame, respectively.\,If \,$K \,=\, I_{H}$\, then the family is called \,$g$-fusion frame and it has been widely studied in \cite{P, Ghosh, G, Ahmadi}.
\end{definition}

Define the space
\[l^{\,2}\left(\,\left\{\,H_{j}\,\right\}_{ j \,\in\, J}\,\right) \,=\, \left \{\,\{\,f_{\,j}\,\}_{j \,\in\, J} \,:\, f_{\,j} \;\in\; H_{j},\; \sum\limits_{\,j \,\in\, J}\, \left \|\,f_{\,j}\,\right \|^{\,2} \,<\, \infty \,\right\}\]
with inner product is given by 
\[\left<\,\{\,f_{\,j}\,\}_{ j \,\in\, J} \,,\, \{\,g_{\,j}\,\}_{ j \,\in\, J}\,\right> \;=\; \sum\limits_{\,j \,\in\, J}\, \left<\,f_{\,j} \,,\, g_{\,j}\,\right>_{H_{j}}.\]\,Clearly \,$l^{\,2}\left(\,\left\{\,H_{j}\,\right\}_{ j \,\in\, J}\,\right)$\; is a Hilbert space with the pointwise operations \cite{Sadri}. 

\subsection{Controlled $K$-$g$-fusion frame}
\smallskip\hspace{.6 cm}
Controlled frame is one of the newest generalization of frame.\,P. Balaz et al.\,\cite{B} introduced controlled frame to improve the numerical efficiency of interactive algorithms for inverting the frame operator.\,In recent times, several generalizations of controlled frame namely, controlled\,$K$-frame \cite{N}, controlled\,$g$-frame \cite{F}, controlled fusion frame \cite{AK}, controlled $g$-fusion frame \cite{HS}, controlled $K$-$g$-fusion frame \cite{GR} etc. have been appeared.

\begin{definition}\cite{GR}
Let \,$K \,\in\, \mathcal{B}\,(\,H\,)$\, and \,$\left\{\,W_{j}\,\right\}_{ j \,\in\, J}$\, be a collection of closed subspaces of \,$H$\, and \,$\left\{\,v_{j}\,\right\}_{ j \,\in\, J}$\, be a collection of positive weights.\,Let \,$\left\{\,H_{j}\,\right\}_{ j \,\in\, J}$\, be a sequence of Hilbert spaces, \,$T,\, U \,\in\, \mathcal{G}\,\mathcal{B}\,(\,H\,)$\, and \,$\Lambda_{j} \,\in\, \mathcal{B}\,(\,H,\, H_{j}\,)$\, for each \,$j \,\in\, J$.\,Then the family \,$\Lambda_{T\,U} \,=\, \left\{\,\left(\,W_{j},\, \Lambda_{j},\, v_{j}\,\right)\,\right\}_{j \,\in\, J}$\, is a \,$(\,T,\,U\,)$-controlled $K$-$g$-fusion frame for \,$H$\, if there exist constants \,$0 \,<\, A \,\leq\, B \,<\, \infty$\, such that
\begin{equation}\label{eqn1.1}
A\,\|\,K^{\,\ast}\,f\,\|^{\,2} \,\leq\, \sum\limits_{\,j \,\in\, J}\, v^{\,2}_{j}\,\left<\,\Lambda_{j}\,P_{\,W_{j}}\,U\,f,\,  \Lambda_{j}\,P_{\,W_{j}}\,T\,f\,\right> \,\leq\, \,B\,\|\,f \,\|^{\,2}
\end{equation}
for all \,$f \,\in\, H$.\,If \,$\Lambda_{T\,U}$\, satisfies only the right inequality of (\ref{eqn1.1}) it is called a \,$(\,T,\,U\,)$-controlled $g$-fusion Bessel sequence in \,$H$.
\end{definition}

Let \,$\Lambda_{T\,U}$\, be a \,$(\,T,\,U\,)$-controlled $g$-fusion Bessel sequence in \,$H$\, with a bound \,$B$.\,The synthesis operator \,$T_{C} \,:\, \mathcal{K}_{\,\Lambda_{j}} \,\to\, H$\, is defined as 
\[T_{C}\,\left(\,\left\{\,v_{\,j}\,\left(\,T^{\,\ast}\,P_{\,W_{j}}\, \Lambda_{j}^{\,\ast}\,\Lambda_{j}\,P_{\,W_{j}}\,U\,\right)^{1 \,/\, 2}\,f\,\right\}_{j \,\in\, J}\,\right) \,=\, \sum\limits_{\,j \,\in\, J}\,v^{\,2}_{j}\,T^{\,\ast}\,P_{\,W_{j}}\, \Lambda_{j}^{\,\ast}\,\Lambda_{j}\,P_{\,W_{j}}\,U\,f,\]for all \,$f \,\in\, H$\, and the analysis operator \,$T^{\,\ast}_{C} \,:\, H \,\to\, \mathcal{K}_{\,\Lambda_{j}}$\,is given by 
\[T_{C}^{\,\ast}\,f \,=\,  \left\{\,v_{\,j}\,\left(\,T^{\,\ast}\,P_{\,W_{j}}\, \Lambda_{j}^{\,\ast}\,\Lambda_{j}\,P_{\,W_{j}}\,U\,\right)^{1 \,/\, 2}\,f\,\right\}_{j \,\in\, J}\; \;\forall\; f \,\in\, H,\]
where 
\[\mathcal{K}_{\,\Lambda_{j}} \,=\, \left\{\,\left\{\,v_{\,j}\,\left(\,T^{\,\ast}\,P_{\,W_{j}}\, \Lambda_{j}^{\,\ast}\,\Lambda_{j}\,P_{\,W_{j}}\,U\,\right)^{1 \,/\, 2}\,f\,\right\}_{j \,\in\, J} \,:\, f \,\in\, H\,\right\} \,\subset\, l^{\,2}\left(\,\left\{\,H_{j}\,\right\}_{ j \,\in\, J}\,\right).\]
The frame operator \,$S_{C} \,:\, H \,\to\, H$\; is defined as follows:
\[S_{C}\,f \,=\, T_{C}\,T_{C}^{\,\ast}\,f \,=\, \sum\limits_{\,j \,\in\, J}\, v_{j}^{\,2}\,T^{\,\ast}\,P_{\,W_{j}}\, \Lambda_{j}^{\,\ast}\,\Lambda_{j}\,P_{\,W_{j}}\,U\,f\; \;\forall\; f \,\in\, H\]and it is easy to verify that 
\[\left<\,S_{C}\,f,\, f\,\right> \,=\, \sum\limits_{\,j \,\in\, J}\, v^{\,2}_{j}\,\left<\,\Lambda_{j}\,P_{\,W_{j}}\,U\,f,\,  \Lambda_{j}\,P_{\,W_{j}}\,T\,f\,\right>\; \;\forall\; f \,\in\, H.\]

Furthermore, if \,$\Lambda_{T\,U}$\, is a \,$(\,T,\,U\,)$-controlled $K$-$g$-fusion frame with bounds \,$A$\, and \,$B$\, then \,$A\,K\,K^{\,\ast} \,\leq\,S_{C} \,\leq\, B\,I_{H}$.

\subsection{$gc$-fusion frame}
\smallskip\hspace{.6 cm}
Frames and their generalizations are mostly considered in the discrete case. However, they also have been studied in continuous case and provided interesting mathematical probelems.\,Continuous frames were proposed by Kaiser \cite{Ka} and it was independently studied by Ali et al.\,\cite{Al}.\,Faroughi et al.\,\cite{MF} introduced the continuous version of $g$-fusion frame.

\begin{definition}\cite{MF}
Let \,$F \,:\, X \,\to\, \mathbb{H}$\; be such that for each \,$h \,\in\, H$, the mapping \,$x \,\to\, P_{\,F\,(\,x\,)}\,(\,h\,)$\; is measurable (\,i.\,e. is weakly measurable\,) and \,$v \,:\, X \,\to\, \mathbb{R}^{\,+}$\, be a measurable function and let \,$\left\{\,K_{x}\,\right\}_{x \,\in\, X}$\, be a collection of Hilbert spaces.\;For each \,$x \,\in\, X$, suppose that \,$\,\Lambda_{x} \,\in\, \mathcal{B}\,(\,F\,(\,x\,) \,,\, K_{x}\,)$.\;Then \,$\Lambda_{F} \,=\, \left\{\,\left(\,F\,(\,x\,),\, \Lambda_{x},\, v\,(\,x\,)\,\right)\,\right\}_{x \,\in\, X}$\; is called a generalized continuous fusion frame or a gc-fusion frame for \,$H$\, with respect to \,$(\,X,\, \mu\,)$\, and \,$v$, if there exists \,$0 \,<\, A \,\leq\, B \,<\, \infty$\; such that
\[A\, \|\,h\,\|^{\,2} \,\leq\, \int\limits_{\,X}\, v^{\,2}\,(\,x\,)\, \left\|\,\Lambda_{x}\,P_{\,F\,(\,x\,)}\,(\,h\,)\,\right\|^{\,2}\,d\mu \,\leq\, B\, \|\,h\,\|^{\,2}\;\; \;\forall\, h \,\in\, H,\]where \,$P_{\,F\,(\,x\,)}$\, is the orthogonal projection onto the subspace \,$F\,(\,x\,)$.\;$\Lambda_{F}$\, is called a tight gc-fusion frame for \,$H$\, if \,$A \,=\, B$\, and Parseval if \,$A \,=\, B \,=\, 1$.\;If we have only the upper bound, we call \,$\Lambda_{F}$\, is a Bessel gc-fusion mapping for \,$H$.
\end{definition}

Let \,$K \,=\, \oplus_{x \,\in\, X}\,K_{x}$\, and \,$L^{\,2}\left(\,X,\, K\,\right)$\, be a collection of all measurable functions \,$\varphi \,:\, X \,\to\, K$\, such that for each \,$x \,\in\, X,\, \varphi\,(\,x\,) \,\in\, K_{x}$\; and \,$\int\limits_{\,X}\,\left\|\,\varphi\,(\,x\,)\,\right\|^{\,2}\,d\mu \,<\, \infty$. It can be verified that \,$L^{\,2}\left(\,X,\, K\,\right)$\; is a Hilbert space with inner product given by
\[\left<\,\phi,\, \varphi\,\right> \,=\, \int\limits_{\,X}\, \left<\,\phi\,(\,x\,),\, \varphi\,(\,x\,)\,\right>\,d\mu\] for \,$\phi,\, \varphi \,\in\, L^{\,2}\left(\,X,\, K\,\right)$.

\subsection{Continuous controlled $g$-fusion frame}
 \smallskip\hspace{.6 cm}In recent times, controlled frames and their generalizations are also studied in continuous case by many researchers.\,P. Ghosh and T. K. Samanta studied continuous version of controlled $g$-fusion frame in \cite{GG}. 
\begin{definition}\cite{GG}
Let \,$F \,:\, X \,\to\, \mathbb{H}$\, be a mapping, \,$v \,:\, X \,\to\, \mathbb{R}^{\,+}$\, be a measurable function and \,$\left\{\,K_{x}\,\right\}_{x \,\in\, X}$\, be a collection of Hilbert spaces.\;For each \,$x \,\in\, X$, suppose that \,$\,\Lambda_{x} \,\in\, \mathcal{B}\,(\,F\,(\,x\,),\, K_{x}\,)$\, and \,$T,\, U \,\in\, \mathcal{G}\,\mathcal{B}^{\,+}\,(\,H\,)$.\,Then \,$\Lambda_{T\,U} \,=\, \left\{\,\left(\,F\,(\,x\,),\, \Lambda_{x},\, v\,(\,x\,)\,\right)\,\right\}_{x \,\in\, X}$\, is called a continuous \,$(\,T,\,U\,)$-controlled generalized fusion frame or continuous \,$(\,T,\,U\,)$-controlled $g$-fusion frame for \,$H$\, with respect to \,$(\,X,\, \mu\,)$\, and \,$v$, if
\begin{description}
\item[$(i)$]for each \,$f \,\in\, H$, the mapping \,$x \,\to\, P_{F\,(\,x\,)}\,(\,f\,)$\; is measurable (\,i.\,e. is weakly measurable\,).
\item[$(ii)$]there exist constants \,$0 \,<\, A \,\leq\, B \,<\, \infty$\, such that
\begin{equation}\label{eqt1}
A\,\|\,f\,\|^{\,2} \,\leq\, \int\limits_{\,X}\,v^{\,2}\,(\,x\,)\,\left<\,\Lambda_{x}\,P_{\,F\,(\,x\,)}\,U\,f,\, \Lambda_{x}\,P_{\,F\,(\,x\,)}\,T\,f\,\right>\,d\mu_{x} \,\leq\, B\,\|\,f\,\|^{\,2},
\end{equation}
for all \,$f \,\in\, H$, where \,$P_{\,F\,(\,x\,)}$\, is the orthogonal projection onto the subspace \,$F\,(\,x\,)$.\,The constants \,$A,\,B$\, are called the frame bounds.\,If only the right inequality of (\ref{eqt1}) holds then \,$\Lambda_{T\,U}$\, is called a continuous \,$(\,T,\,U\,)$-controlled $g$-fusion Bessel family for \,$H$.
\end{description}
\end{definition}
 
Let \,$\Lambda_{T\,U}$\, be a continuous \,$(\,T,\,U\,)$-controlled $g$-fusion Bessel family for \,$H$.\,Then the operator \,$S_{C} \,:\, H \,\to\,\ H$\, defined by
\[\left<\,S_{C}\,f,\, g\,\right> \,=\, \int\limits_{\,X}\,v^{\,2}\,(\,x\,)\,\left<\,T^{\,\ast}\,P_{F\,(\,x\,)}\,\Lambda_{x}^{\,\ast}\,\Lambda_{x}\,P_{F\,(\,x\,)}\,U\,f,\, g\,\right>\,d\mu_{x}\; \;\forall\; f,\, g \,\in\, H.\]
is called the frame operator.\,If \,$\Lambda_{T\,U}$\, is a continuous \,$(\,T,\,U\,)$-controlled $g$-fusion frame for \,$H$\, then from (\ref{eqt1}), we get
\[A\,\left<\,f,\, f\,\right> \,\leq\, \left<\,S_{C}\,f,\, f\,\right> \,\leq\, B\,\left<\,f,\, f\,\right>\; \;\forall\; f \,\in\, H.\]The bounded linear operator \,$T_{C} \,:\, L^{\,2}\left(\,X,\, K\,\right) \,\to\, H$\, defined by 
\[\left<\,T_{C}\,\Phi,\, g\,\right> \,=\, \int\limits_{\,X}\,v^{\,2}\,(\,x\,)\,\left<\,T^{\,\ast}\,P_{F\,(\,x\,)}\,\Lambda_{x}^{\,\ast}\,\Lambda_{x}\,P_{F\,(\,x\,)}\,U\,f,\, g\,\right>\,d\mu_{x},\]
where for all \,$f \,\in\, H$, \,$\Phi \,=\, \left\{\,v\,(\,x\,)\,\left(\,T^{\,\ast}\,P_{F\,(\,x\,)}\,\Lambda_{x}^{\,\ast}\,\Lambda_{x}\,P_{F\,(\,x\,)}\,U\,\right)^{1 \,/\, 2}\,f\,\right\}_{x \,\in\, X}$\, and \,$g \,\in\, H$, is called synthesis operator and its adjoint operator is called analysis operator.

\subsection{Weaving frame}
\smallskip\hspace{.6 cm}
Woven frame is a new notion in frame theory which has been introduced by Bemrose et al.\,\cite{BM}.\,Two frames \,$\left\{\,f_{\,i}\,\right\}_{i \,\in\, I}$\, and \,$\left\{\,g_{\,i}\,\right\}_{i \,\in\, I}$\, for \,$H$\, are called woven if there exist constants \,$0 \,<\, A \,\leq\, B \,<\, \infty$\, such that for any subset \,$\sigma \,\subset\, I$\, the family \,$\left\{\,f_{\,i}\,\right\}_{i \,\in\, \sigma} \,\cup\, \left\{\,g_{\,i}\,\right\}_{i \,\in\, \sigma^{\,c}}$\, is a frame for \,$H$.\,This frame has been generalized for the discrete as well as the continuous case such as woven fusion frame \cite{SG}, woven \,$g$-frame \cite{DL}, woven \,$g$-fusion frame \cite{MM}, woven \,$K$-$g$-fusion frame \cite{VG}, Continuous weaving frame \cite{LK}, Continuous weaving fusion frame \cite{RA} etc. \\

In this paper,\, woven continuous controlled \,$K$-$g$-fusion frame in Hilbert spaces is presented and some of their properties are going to be established.\,We discuss sufficient conditions for weaving continuous controlled \,$K$-$g$-fusion frame.\,Construction of woven continuous controlled \,$K$-$g$-fusion frame by bounded linear operator is given.\,At the end, we discuss a perturbation result of woven continuous controlled $K$-$g$-fusion frame.

\section{Weaving continuous controlled $K$-$g$-fusion frame }

\smallskip\hspace{.6 cm}In this section, we first give the continuous version of controlled \,$K$-$g$-fusion frame for \,$H$\, and then present weaving continuous controlled \,$K$-$g$-fusion frame for \,$H$. 

\begin{definition}
Let \,$K \,\in\, \mathcal{B}\,(\,H\,)$\, and \,$F \,:\, X \,\to\, \mathbb{H}$\, be a mapping, \,$v \,:\, X \,\to\, \mathbb{R}^{\,+}$\, be a measurable function and \,$\left\{\,K_{x}\,\right\}_{x \,\in\, X}$\, be a collection of Hilbert spaces.\;For each \,$x \,\in\, X$, suppose that \,$\,\Lambda\,(\,x\,) \,\in\, \mathcal{B}\,(\,F\,(\,x\,),\, K_{x}\,)$\, and \,$T,\, U \,\in\, \mathcal{G}\,\mathcal{B}^{\,+}\,(\,H\,)$.\,Then \,$\Lambda_{T\,U} \,=\, \left\{\,\left(\,F\,(\,x\,),\, \Lambda\,(\,x\,),\, v\,(\,x\,)\,\right)\,\right\}_{x \,\in\, X}$\, is called a continuous \,$(\,T,\,U\,)$-controlled $K$-$g$-fusion frame for \,$H$\, with respect to \,$(\,X,\, \mu\,)$\, and \,$v$, if
\begin{description}
\item[$(i)$]for each \,$f \,\in\, H$, the mapping \,$x \,\to\, P_{F\,(\,x\,)}\,(\,f\,)$\; is measurable (\,i.\,e. is weakly measurable\,).
\item[$(ii)$]there exist constants \,$0 \,<\, A \,\leq\, B \,<\, \infty$\, such that
\begin{equation}\label{eq1}
A\,\left\|\,K^{\,\ast}\,f\,\right\|^{\,2} \,\leq\, \int\limits_{\,X}\,v^{\,2}\,(\,x\,)\,\left<\,\Lambda\,(\,x\,)\,P_{\,F\,(\,x\,)}\,U\,f,\, \Lambda\,(\,x\,)\,P_{\,F\,(\,x\,)}\,T\,f\,\right>\,d\mu_{x} \,\leq\, B\,\|\,f\,\|^{\,2},
\end{equation}
for all \,$f \,\in\, H$, where \,$P_{\,F\,(\,x\,)}$\, is the orthogonal projection onto the subspace \,$F\,(\,x\,)$.\,The constants \,$A,\,B$\, are called the frame bounds.
\end{description}
\end{definition}
Now, we consider the following cases:
\begin{description}
\item[$(I)$]If only the right inequality of (\ref{eq1}) holds then \,$\Lambda_{T\,U}$\, is called a continuous \,$(\,T,\,U\,)$-controlled $K$-$g$-fusion Bessel family for \,$H$.
\item[$(II)$]If \,$U \,=\, I_{H}$\, then \,$\Lambda_{T\,U}$\, is called a continuous \,$(\,T,\,I_{H}\,)$-controlled $K$-$g$-fusion frame for \,$H$. 
\item[$(III)$]If \,$T \,=\, U \,=\, I_{H}$\, then \,$\Lambda_{T\,U}$\, is called a continuous $K$-$g$-fusion frame for \,$H$.
\item[$(IV)$]if \,$K \,=\, I_{H}$\, then \,$\Lambda_{T\,U}$\, is called a continuous \,$(\,T,\,U\,)$-controlled $g$-fusion frame for \,$H$.  
\end{description}
 
\begin{remark}
If the measure space \,$X \,=\, \mathbb{N}$\, and \,$\mu$\, is the counting measure then a continuous \,$(\,T,\,U\,)$-controlled $K$-$g$-fusion frame will be the discrete \,$(\,T,\,U\,)$-controlled $K$-$g$-fusion frame.   
\end{remark} 

\subsubsection{Example}
Let \,$H \,=\, \mathbb{R}^{\,3}$\, and \,$\left\{\,e_{\,1},\,e_{\,2},\, e_{\,3}\,\right\}$\, be an standard orthonormal basis for \,$H$.\,Consider
\[\mathcal{B} \,=\, \left\{\,x \,\in\, \mathbb{R}^{\,3} \,:\, \|\,x\,\| \,\leq\, 1\,\right\}.\]
Then it is a measure space equipped with the Lebesgue measure \,$\mu$.\,Suppose \,$\left\{\,B_{\,1},\,B_{\,2},\, B_{\,3}\,\right\}$\, is a partition of \,$\mathcal{B}$\, where \,$\mu\,(\,B_{1}\,) \,\geq\, \mu\,(\,B_{2}\,) \,\geq\, \mu\,(\,B_{3}\,) \,>\, 1$.\,Let \,$\mathbb{H} \,=\, \left\{\,W_{\,1},\,W_{\,2},\, W_{\,3}\,\right\}$, where \,$W_{1} \,=\, \overline{span}\,\left\{\,e_{\,2},\, e_{\,3}\,\right\}$, \,$W_{2} \,=\, \overline{span}\,\left\{\,e_{\,1},\, e_{\,3}\,\right\}$\, and \,$W_{3} \,=\, \overline{span}\,\left\{\,e_{\,1},\, e_{\,2}\,\right\}$.\,Define 
\[F \,:\, \mathcal{B} \,\to\, \mathbb{H}\hspace{.5cm}\text{by} \hspace{.3cm} F\,(\,x\,) \,=\, \begin{cases}
W_{1} & \text{if\;\;}\; x \,\in\, B_{1} \\ W_{2} & \text{if\;\;}\; x \,\in\, B_{2}\\ W_{3} & \text{if\;\;}\; x \,\in\, B_{3} \end{cases}\] 
and
\[v \,:\, \mathcal{B} \,\to\, [\,0\, \infty)\hspace{.5cm}\text{by} \hspace{.3cm} v\,(\,x\,) \,=\, \begin{cases}
1 & \text{if\;\;}\; x \,\in\, B_{1} \\ 2 & \text{if\;\;}\; x \,\in\, B_{2}\\ \,-\, 1 & \text{if\;\;}\; x \,\in\, B_{3}\,. \end{cases}\]
It is easy to verify that \,$F$\, and \,$v$\, are measurable functions.\,For each \,$x \,\in\, \mathcal{B}$, define the operators
\[\Lambda\,(\,x\,)\,(\,f\,) \,=\, \dfrac{1}{\sqrt{\mu\,(\,B_{k}\,)}}\left<\,f,\, e_{k}\,\right>\,e_{k},\; \;f \,\in\, H,\]where \,$k$\, is such that \,$x \,\in\, \mathcal{B}_{k}$\, and \,$K \,:\, H \,\to\, H$\, by
\[K\,e_{\,1} \,=\, e_{\,1},\, K\,e_{\,2} \,=\, e_{\,2},\, K\,e_{\,3} \,=\, 0.\]
It is easy to verify that \,$K^{\,\ast}\,e_{\,1} \,=\, e_{\,1},\, K^{\,\ast}\,e_{\,2} \,=\, e_{\,2},\, K^{\,\ast}\,e_{\,3} \,=\, 0$.\,Now, for any \,$f \,\in\, H$, we have
\[\left\|\,K^{\,\ast}\,f\,\right\|^{\,2} \,=\, \left\|\,\sum\limits_{i \,=\, 1}^{\,3}\,\left<\,f,\, e_{k}\,\right>\,K^{\,\ast}\,e_{k}\,\right\|^{\,2} \,=\, \left|\,\left<\,f,\, e_{1}\,\right>\,\right|^{\,2} \,+\, \left|\,\left<\,f,\, e_{2}\,\right>\,\right|^{\,2} \,\leq\, \|\,f\,\|^{\,2}.\]
Let \,$T\,\left(\,f_{\,1},\, f_{\,2},\, f_{\,3}\,\right) \,=\, \left(\,5\,f_{\,1},\, 4\,f_{\,2},\, 5\,f_{\,3}\,\right)$\, and \,$U\,\left(\,f_{\,1},\, f_{\,2},\, f_{\,3}\,\right) \,=\, \left(\,\dfrac{f_{\,1}}{6},\, \dfrac{f_{\,2}}{3},\, \dfrac{f_{\,3}}{6}\,\right)$\, be two operators on \,$H$.\,Then it is easy to verify that \,$T,\, U \,\in\, \mathcal{G}\,\mathcal{B}^{\,+}\,(\,H\,)$\, and $\,T\,U \,=\, U\,T$.\,Now, for any \,$f \,=\, \left(\,f_{\,1},\, f_{\,2},\, f_{\,3}\,\right) \,\in\, H$, we have 
\begin{align*}
&\int\limits_{\,\mathcal{B}}\,v^{\,2}\,(\,x\,)\,\left<\,\Lambda\,(\,x\,)\,P_{\,F\,(\,x\,)}\,U\,f,\, \Lambda\,(\,x\,)\,P_{\,F\,(\,x\,)}\,T\,f\,\right>\,d\mu_{x}\\
&=\,\sum\limits_{i \,=\, 1}^{\,3}\int\limits_{\,\mathcal{B}_{i}}\,v^{\,2}\,(\,x\,)\,\left<\,\Lambda\,(\,x\,)\,P_{\,F\,(\,x\,)}\,U\,f,\, \Lambda\,(\,x\,)\,P_{\,F\,(\,x\,)}\,T\,f\,\right>\,d\mu_{x}\\
&=\,\dfrac{5}{6}f_{\,1}^{\,2} \,+\, \dfrac{16}{3}\,f_{\,2}^{\,2} \,+\, \dfrac{5}{6}\,f_{\,3}^{\,2}.\\
&\Rightarrow\,\left\|\,K^{\,\ast}\,f\,\right\|^{\,2} \,\leq\, \int\limits_{\,\mathcal{B}}\,v^{\,2}\,(\,x\,)\,\left<\,\Lambda\,(\,x\,)\,P_{\,F\,(\,x\,)}\,U\,f,\, \Lambda\,(\,x\,)\,P_{\,F\,(\,x\,)}\,T\,f\,\right>\,d\mu_{x} \,\leq\, \dfrac{16}{3}\,\|\,f\,\|^{\,2}. 
\end{align*}
Thus, \,$\Lambda_{T\,U}$\, be a continuous \,$(\,T,\,U\,)$-controlled $K$-$g$-fusion frame for \,$\mathbb{R}^{\,3}$.\\

Now, we present woven continuous controlled $K$-$g$-fusion frame for \,$H$.
 
\begin{definition}
A family of continuous \,$(\,T,\,U\,)$-controlled $K$-$g$-fusion frames \,$\left\{\,\left(\,F_{\,i}\,(\,x\,),\, \Lambda_{i}\,(\,x\,),\, v_{\,i}\,(\,x\,)\,\right)\,\right\}_{i \,\in\, [\,m\,],\, x \,\in\, X}$\, for the Hilbert space \,$H$\, is said to be woven continuous \,$(\,T,\,U\,)$-controlled $K$-$g$-fusion frame if there exist universal positive constants \,$0 \,<\, A \,\leq\, B \,<\, \infty$\, such that for each partition \,$\left\{\,\sigma_{\,i}\,\right\}_{i \,\in\, [\,m\,]}$\, of \,$X$, the family \,$\left\{\,\left(\,F_{\,i}\,(\,x\,),\, \Lambda_{i}\,(\,x\,),\, v_{\,i}\,(\,x\,)\,\right)\,\right\}_{i \,\in\, [\,m\,],\, x \,\in\, \sigma_{\,i}}$\, is a continuous \,$(\,T,\,U\,)$-controlled $K$-$g$-fusion frame for \,$H$\, with bounds \,$A$\, and \,$B$.
\end{definition}

Each family \,$\left\{\,\left(\,F_{\,i}\,(\,x\,),\, \Lambda_{i}\,(\,x\,),\, v_{\,i}\,(\,x\,)\,\right)\,\right\}_{i \,\in\, [\,m\,],\, x \,\in\, \sigma_{\,i}}$\, is called a weaving continuous \,$(\,T,\,U\,)$-controlled $K$-$g$-fusion frame.\,For abbreviation, we use W. C. C. K. G. F. F. instead of the statement of woven continuous \,$(\,T,\,U\,)$-controlled $K$-$g$-fusion frame.          

In the following proposition, we will see that every woven continuous controlled \,$K$-$g$-fusion frame has a universal upper bound.

\begin{proposition}
Let \,$\left\{\,\left(\,F_{\,i}\,(\,x\,),\, \Lambda_{i}\,(\,x\,),\, v_{\,i}\,(\,x\,)\,\right)\,\right\}_{x \,\in\, X}$\, be a continuous \,$(\,T,\,U\,)$-controlled $K$-$g$-fusion Bessel family for \,$H$\, with bound \,$B_{\,i}$\, for each \,$i \,\in\, [\,m\,]$.\,Then for any partition \,$\left\{\,\sigma_{\,i}\,\right\}_{i \,\in\, [\,m\,]}$\, of \,$X$, the family \,$\left\{\,\left(\,F_{\,i}\,(\,x\,),\, \Lambda_{i}\,(\,x\,),\, v_{\,i}\,(\,x\,)\,\right)\,\right\}_{i \,\in\, [\,m\,],\, x \,\in\, \sigma_{\,i}}$\, is a continuous \,$(\,T,\,U\,)$-controlled $K$-$g$-fusion Bessel family for \,$H$.   
\end{proposition}

\begin{proof}
Let \,$\left\{\,\sigma_{\,i}\,\right\}_{i \,\in\, [\,m\,]}$\, be a arbitrary partition of \,$X$.\,For each \,$f \,\in\, H$, we have
\begin{align*}
&\sum\limits_{i \,\in\, [\,m\,]}\,\int\limits_{\,\sigma_{\,i}}\,v_{\,i}^{\,2}\,(\,x\,)\,\left<\,\Lambda_{i}\,(\,x\,)\,P_{\,F_{\,i}\,(\,x\,)}\,U\,f,\, \Lambda_{i}\,(\,x\,)\,P_{\,F_{\,i}\,(\,x\,)}\,T\,f\,\right>\,d\mu_{x}\\
&\leq\,\sum\limits_{i \,\in\, [\,m\,]}\,\int\limits_{\,X}\,v_{\,i}^{\,2}\,(\,x\,)\,\left<\,\Lambda_{i}\,(\,x\,)\,P_{\,F_{\,i}\,(\,x\,)}\,U\,f,\, \Lambda_{i}\,(\,x\,)\,P_{\,F_{\,i}\,(\,x\,)}\,T\,f\,\right>\,d\mu_{x}\\
&\leq\,\left(\,\sum\limits_{i \,\in\, [\,m\,]}\,B_{\,i}\,\right)\,\|\,f\,\|^{\,2}. 
\end{align*}
This completes the proof. \\
\end{proof}
Next, we gives a characterization of W. C. C. K. G. F. F. for \,$H$\, in terms of an operator.
 
\begin{theorem}
Let the families \,$\Lambda \,=\, \left\{\,\left(\,F\,(\,x\,),\, \Lambda\,(\,x\,),\, v\,(\,x\,)\,\right)\,\right\}_{x \,\in\, X}$\,  and \,$\Gamma \,=\, \left\{\,\left(\,G\,(\,x\,),\, \Lambda\,(\,x\,),\, v\,(\,x\,)\,\right)\,\right\}_{x \,\in\, X}$\, be continuous \,$(\,T,\,U\,)$-controlled \,$K$-$g$-fusion frames for \,$H$.The the following statements are equivalent:
\begin{itemize}
\item[$(I)$]$\Lambda$\, and \,$\Gamma$ are W. C. C. K. G. F. F. for \,$H$. 
\item[$(II)$]For each partition \,$\sigma$\, of \,$X$, there exist \,$\alpha \,>\, 0$\, and a bounded linear operator \,$\Theta_{\sigma} \,:\, L^{\,2}_{\,\sigma}\left(\,X,\, K\,\right) \,\to\, H$\, defined by
\begin{align*}
\left<\,\Theta_{\sigma}\,\Phi,\, g\,\right>& \,=\, \int\limits_{\,\sigma}\,v^{\,2}\,(\,x\,)\,\left<\,T^{\,\ast}\,P_{F\,(\,x\,)}\,\Lambda\,(\,x\,)^{\,\ast}\,\Lambda\,(\,x\,)\,P_{F\,(\,x\,)}\,U\,f,\, g\,\right>\,d\mu_{x} \,+\\
&\int\limits_{\,\sigma^{\,c}}\,v^{\,2}\,(\,x\,)\,\left<\,T^{\,\ast}\,P_{G\,(\,x\,)}\,\Gamma\,(\,x\,)^{\,\ast}\,\Gamma\,(\,x\,)\,P_{G\,(\,x\,)}\,U\,f,\, g\,\right>\,d\mu_{x},\, \,g \,\in\, H
\end{align*} 
such that \,$\alpha\,K\,K^{\,\ast} \,\leq\, \Theta_{\sigma}\,\Theta_{\sigma}^{\,\ast}$, where
\[L^{\,2}_{\,\sigma}\left(\,X,\, K\,\right) \,=\, \left\{\,\Phi \,=\, \phi \,\cup\, \psi \,:\, \int\limits_{\,X}\,\|\,\Phi\,\|^{\,2}\,d\,\mu \,<\, \infty\,\right\},\]
where for all \,$f \,\in\, H$, \,$\phi \,=\, \left\{\,v\,(\,x\,)\,\left(\,T^{\,\ast}\,P_{F\,(\,x\,)}\,\Lambda\,(\,x\,)^{\,\ast}\,\Lambda\,(\,x\,)\,P_{F\,(\,x\,)}\,U\,\right)^{1 \,/\, 2}\,f\,\right\}_{x \,\in\, \sigma}$\, and \,$\psi \,=\, \left\{\,v\,(\,x\,)\,\left(\,T^{\,\ast}\,P_{G\,(\,x\,)}\,\Gamma\,(\,x\,)^{\,\ast}\,\Gamma\,(\,x\,)\,P_{G\,(\,x\,)}\,U\,\right)^{1 \,/\, 2}\,f\,\right\}_{x \,\in\, \sigma^{\,c}}$.
\end{itemize} 
\end{theorem}

\begin{proof}$(I) \,\Rightarrow\, (II)$
Suppose that \,$A$\, and \,$B$\, are the universal lower and upper bounds for $\Lambda$\, and \,$\Gamma$.\,Take \,$\Theta_{\sigma} \,=\, T_{C}^{\,\sigma}$, for every partition \,$\sigma$\, of \,$X$, where \,$T_{C}^{\,\sigma}$\, is the synthesis operator of 
\[\left\{\,\left(\,F\,(\,x\,),\, \Lambda\,(\,x\,),\, v\,(\,x\,)\,\right)\,\right\}_{x \,\in\, \sigma} \,\cup\, \left\{\,\left(\,G\,(\,x\,),\, \Lambda\,(\,x\,),\, v\,(\,x\,)\,\right)\,\right\}_{x \,\in\, \sigma^{\,c}}.\]
Thus, for each \,$\Phi \,\in\, L^{\,2}_{\,\sigma}\left(\,X,\, K\,\right)$, we have
\begin{align*}
\left<\,\Theta_{\sigma}\,\Phi,\, g\,\right> &\,=\, \left<\,T_{C}^{\,\sigma}\,\Phi,\, g\,\right> \,=\, \int\limits_{\,\sigma}\,v^{\,2}\,(\,x\,)\,\left<\,T^{\,\ast}\,P_{F\,(\,x\,)}\,\Lambda\,(\,x\,)^{\,\ast}\,\Lambda\,(\,x\,)\,P_{F\,(\,x\,)}\,U\,f,\, g\,\right>\,d\mu_{x} \,+\\
&\int\limits_{\,\sigma^{\,c}}\,v^{\,2}\,(\,x\,)\,\left<\,T^{\,\ast}\,P_{G\,(\,x\,)}\,\Gamma\,(\,x\,)^{\,\ast}\,\Gamma\,(\,x\,)\,P_{G\,(\,x\,)}\,U\,f,\, g\,\right>\,d\mu_{x},\, \,g \,\in\, H. 
\end{align*}
Since \,$\Lambda$\, and \,$\Gamma$ are woven, for each \,$f \,\in\, H$, we have
\begin{align*}
&A\,\left\|\,K^{\,\ast}\,f\,\right\|^{\,2} \,\leq\, \left\|\,\left(\,T_{C}^{\,\sigma}\,\right)^{\,\ast}\,f\,\right\|^{\,2} \,=\, \left\|\,\Theta_{\sigma}^{\,\ast}\,f\,\right\|^{\,2}\\
&\Rightarrow\,\alpha\,K\,K^{\,\ast} \,\leq\, \Theta_{\sigma}\,\Theta_{\sigma}^{\,\ast},\, \,\alpha \,=\, A.  
\end{align*} 
$(II) \,\Rightarrow\, (I)$\, Let \,$\sigma$\, be a partition of \,$X$\, and \,$f \,\in\, H$.\,Now it is easy to verify that
\begin{align*}
\Theta_{\sigma}^{\,\ast}\,f \,=\, &\left\{\,v\,(\,x\,)\,\left(\,T^{\,\ast}\,P_{F\,(\,x\,)}\,\Lambda\,(\,x\,)^{\,\ast}\,\Lambda\,(\,x\,)\,P_{F\,(\,x\,)}\,U\,\right)^{1 \,/\, 2}\,f\,\right\}_{x \,\in\, \sigma} \,\cup\, \\
&\left\{\,v\,(\,x\,)\,\left(\,T^{\,\ast}\,P_{G\,(\,x\,)}\,\Gamma\,(\,x\,)^{\,\ast}\,\Gamma\,(\,x\,)\,P_{G\,(\,x\,)}\,U\,\right)^{1 \,/\, 2}\,f\,\right\}_{x \,\in\, \sigma^{\,c}}
\end{align*}
Thus, for each \,$f \,\in\, H$, we have 
\begin{align*}
\alpha\,\left\|\,K^{\,\ast}\,f\,\right\|^{\,2} \,\leq\, \left\|\,\Theta_{\sigma}^{\,\ast}\,f\,\right\|^{\,2} \,=\,& \int\limits_{\,\sigma}\,v^{\,2}\,(\,x\,)\,\left<\,\Lambda\,(\,x\,)\,P_{\,F\,(\,x\,)}\,U\,f,\, \Lambda\,(\,x\,)\,P_{\,F\,(\,x\,)}\,T\,f\,\right>\,d\mu_{x} \,+\\
&\int\limits_{\,\sigma^{\,c}}\,v^{\,2}\,(\,x\,)\,\left<\,\Gamma\,(\,x\,)\,P_{\,G\,(\,x\,)}\,U\,f,\, \Gamma\,(\,x\,)\,P_{\,G\,(\,x\,)}\,T\,f\,\right>\,d\mu_{x}. 
\end{align*}
Hence, $\Lambda$\, and \,$\Gamma$ are W. C. C. K. G. F. F. for \,$H$.
\end{proof}

In the following theorem, we will construct W. C. C. K. G. F. F. for \,$H$\, by using bounded linear operator.   

\begin{theorem}\label{t3.6}
Let \,$\left\{\,\left(\,F_{\,i}\,(\,x\,),\, \Lambda_{i}\,(\,x\,),\, v_{\,i}\,(\,x\,)\,\right)\,\right\}_{i \,\in\, [\,m\,],\, x \,\in\, \sigma_{\,i}}$\, be a W. C. C. K. G. F. F. for \,$H$\, with universal bounds \,$A$\, and \,$B$.\,If \,$V \,\in\, \mathcal{B}\,(\,H\,)$\, be invertible operator such that \,$V^{\,\ast}$\, commutes with \,$T,\,U$\, and \,$V$\, commutes with \,$K$, then the family \,$\left\{\,\left(\,V\,F_{\,i}\,(\,x\,),\, \Lambda_{i}\,(\,x\,)\,P_{F_{\,i}\,(\,x\,)}\,V^{\,\ast},\, v_{\,i}\,(\,x\,)\,\right)\,\right\}_{i \,\in\, [\,m\,],\, x \,\in\, \sigma_{\,i}}$\, is a W. C. C. K. G. F. F. for \,$H$.  
\end{theorem}

\begin{proof}
Since \,$P_{F_{\,i}\,(\,x\,)}\,V^{\,\ast} \,=\, P_{F_{\,i}\,(\,x\,)}\,V^{\,\ast}\,P_{V\,F_{\,i}\,(\,x\,)}$\, for all \,$x \,\in\, \sigma_{\,i}$\, and \,$i \,\in\, [\,m\,]$, the mapping \,$x \,\to\, P_{V\,F_{\,i}\,(\,x\,)}$\, is weakly measurable.\,For each \,$f \,\in\, H$, we have
\begin{align*}
&\sum\limits_{i \,\in\, [\,m\,]}\,\int\limits_{\,X}\,v_{\,i}^{\,2}\,(\,x\,)\,\left<\,\Lambda_{i}\,(\,x\,)\,P_{F_{\,i}\,(\,x\,)}\,V^{\,\ast}\,P_{V\,F_{\,i}\,(\,x\,)}\,U\,f,\, \Lambda_{i}\,(\,x\,)\,P_{F_{\,i}\,(\,x\,)}\,V^{\,\ast}\,P_{V\,F_{\,i}\,(\,x\,)}\,T\,f\,\right>\,d\mu_{x}\\
&=\,\sum\limits_{i \,\in\, [\,m\,]}\,\int\limits_{\,X}\,v_{\,i}^{\,2}\,(\,x\,)\,\left<\,\Lambda_{i}\,(\,x\,)\,P_{F_{\,i}\,(\,x\,)}\,V^{\,\ast}\,U\,f,\, \Lambda_{i}\,(\,x\,)\,P_{F_{\,i}\,(\,x\,)}\,V^{\,\ast}\,T\,f\,\right>\,d\mu_{x}\\
&=\,\sum\limits_{i \,\in\, [\,m\,]}\,\int\limits_{\,X}\,v_{\,i}^{\,2}\,(\,x\,)\,\left<\,\Lambda_{i}\,(\,x\,)\,P_{F_{\,i}\,(\,x\,)}\,U\,V^{\,\ast}\,f,\, \Lambda_{i}\,(\,x\,)\,P_{F_{\,i}\,(\,x\,)}\,T\,V^{\,\ast}\,f\,\right>\,d\mu_{x}\\
&\leq\,B\,\left\|\,V^{\,\ast}\,f\,\right\|^{\,2} \,\leq\, B\,\left\|\,V\,\right\|^{\,2}\,\|\,f\,\|^{\,2}. 
\end{align*}
On the other hand, for each \,$f \,\in\, H$, we have
\begin{align*}
&\sum\limits_{i \,\in\, [\,m\,]}\,\int\limits_{\,X}\,v_{\,i}^{\,2}\,(\,x\,)\,\left<\,\Lambda_{i}\,(\,x\,)\,P_{F_{\,i}\,(\,x\,)}\,V^{\,\ast}\,P_{V\,F_{\,i}\,(\,x\,)}\,U\,f,\, \Lambda_{i}\,(\,x\,)\,P_{F_{\,i}\,(\,x\,)}\,V^{\,\ast}\,P_{V\,F_{\,i}\,(\,x\,)}\,T\,f\,\right>\,d\mu_{x}\\
&\geq\, A\,\left\|\,K^{\,\ast}\,V^{\,\ast}\,f\,\right\|^{\,2} \,=\, A\,\left\|\,V^{\,\ast}\,K^{\,\ast}\,f\,\right\|^{\,2} \,\geq\, A\,\left\|\,V^{\,-\, 1}\,\right\|^{\,-\, 2}\,\left\|\,K^{\,\ast}\,f\,\right\|^{\,2}.
\end{align*} 
This completes the proof.
\end{proof}

\begin{corollary}
Let \,$\left\{\,\left(\,F_{\,i}\,(\,x\,),\, \Lambda_{i}\,(\,x\,),\, v_{\,i}\,(\,x\,)\,\right)\,\right\}_{i \,\in\, [\,m\,],\, x \,\in\, \sigma_{\,i}}$\, be a W. C. C. K. G. F. F. for \,$H$\, with universal bounds \,$A$\, and \,$B$.\,If \,$V \,\in\, \mathcal{B}\,(\,H\,)$\, be invertible operator such that \,$V^{\,\ast}$\, commutes with \,$T,\,U$\, and \,$V$\, commutes with \,$K$, then the family \,$\left\{\,\left(\,V\,F_{\,i}\,(\,x\,),\, \Lambda_{i}\,(\,x\,)\,P_{F_{\,i}\,(\,x\,)}\,V^{\,\ast},\, v_{\,i}\,(\,x\,)\,\right)\,\right\}_{i \,\in\, [\,m\,],\, x \,\in\, \sigma_{\,i}}$\, is a W. C. C. $V$K$V^{\,\ast}$. G. F. F. for \,$H$.  
\end{corollary}

\begin{proof}
According to the proof of theorem \ref{t3.6}, universal upper bounds is \,$B\,\|\,V\,\|^{\,2}$. On the other hand, for each \,$f \,\in\, H$, we have
\begin{align*}
&\dfrac{A}{\|\,V\,\|^{\,2}}\, \left\|\,\left(\,V\,K\,V^{\,\ast}\,\right)^{\,\ast}\,f\,\right\|^{\,2} \,=\, \dfrac{A}{\|\,V\,\|^{\,2}}\, \left\|\,V\,K^{\,\ast}\,V^{\,\ast}\,f\,\right\|^{\,2} \,\leq\, A\, \left\|\,K^{\,\ast}\,V^{\,\ast}\,f\,\right\|^{\,2}\\
&\leq\,\sum\limits_{i \,\in\, [\,m\,]}\,\int\limits_{\,X}\,v_{\,i}^{\,2}\,(\,x\,)\,\left<\,\Lambda_{i}\,(\,x\,)\,P_{F_{\,i}\,(\,x\,)}\,U\,V^{\,\ast}\,f,\, \Lambda_{i}\,(\,x\,)\,P_{F_{\,i}\,(\,x\,)}\,T\,V^{\,\ast}\,f\,\right>\,d\mu_{x}\\
&=\,\sum\limits_{i \,\in\, [\,m\,]}\,\int\limits_{\,X}\,v_{\,i}^{\,2}\,(\,x\,)\,\left<\,\Lambda_{i}\,(\,x\,)\,P_{F_{\,i}\,(\,x\,)}\,V^{\,\ast}\,P_{V\,F_{\,i}\,(\,x\,)}\,U\,f,\, \Lambda_{i}\,(\,x\,)\,P_{F_{\,i}\,(\,x\,)}\,V^{\,\ast}\,P_{V\,F_{\,i}\,(\,x\,)}\,T\,f\,\right>\,d\mu_{x}.
\end{align*}
This completes the proof.
\end{proof}

\begin{theorem}
Let \,$V \,\in\, \mathcal{B}\,(\,H\,)$\, be invertible operator such that \,$V^{\,\ast}, \left(\,V^{\,-\, 1}\,\right)^{\,\ast}$\, commutes with \,$T$\, and \,$U$.\,Suppose \,$\left\{\,\left(\,V\,F_{\,i}\,(\,x\,),\, \Lambda_{i}\,(\,x\,)\,P_{F_{\,i}\,(\,x\,)}\,V^{\,\ast},\, v_{\,i}\,(\,x\,)\,\right)\,\right\}_{i \,\in\, [\,m\,],\, x \,\in\, \sigma_{\,i}}$\, is a W. C. C. K. G. F. F. for \,$H$\, with universal bounds \,$A$\, and \,$B$.\,Then the family \,$\left\{\,\left(\,F_{\,i}\,(\,x\,),\, \Lambda_{i}\,(\,x\,),\, v_{\,i}\,(\,x\,)\,\right)\,\right\}_{i \,\in\, [\,m\,],\, x \,\in\, \sigma_{\,i}}$\, be a W. C. C. $V^{\,-\, 1}$K$V$. G. F. F. for \,$H$.    
\end{theorem}

\begin{proof}
Now, for each \,$f \,\in\, H$, using Theorem \ref{th1.01}, we have
\begin{align*}
&\dfrac{A}{\|\,V\,\|^{\,2}}\,\left \|\,\left(\,V^{\,-\, 1}\,K\,V\,\right)^{\,\ast}\,f\,\right \|^{\,2} \,=\, \dfrac{A}{\|\,V\,\|^{\,2}}\,\left\|\,V^{\,\ast}\,K^{\,\ast}\,(\,V^{\,-\, 1}\,)^{\,\ast}\,f\,\right\|^{\,2}\\
& \,\leq\, A\;\left\|\,K^{\,\ast}\,\left(\,V^{\,-\, 1}\,\right)^{\,\ast}\,f\,\right\|^{\,2}\\
&\leq\sum\limits_{i \,\in\, [\,m\,]}\int\limits_{\,X}v_{i}^{\,2}\,(\,x\,)\left<\,\Lambda_{i}\,(\,x\,)P_{F_{i}\,(\,x\,)}V^{\,\ast}P_{V\,F_{i}\,(\,x\,)}U\left(V^{\,-\, 1}\right)^{\,\ast}f,\, \Lambda_{i}\,(\,x\,)P_{F_{i}\,(\,x\,)}V^{\,\ast}P_{V\,F_{i}\,(\,x\,)}T\left(V^{\,-\, 1}\right)^{\,\ast}f\,\right>\,d\mu_{x}\\
&\leq\,\sum\limits_{i \,\in\, [\,m\,]}\,\int\limits_{\,X}\,v_{\,i}^{\,2}\,(\,x\,)\,\left<\,\Lambda_{i}\,(\,x\,)\,P_{F_{\,i}\,(\,x\,)}\,V^{\,\ast}\,U\,\left(\,V^{\,-\, 1}\,\right)^{\,\ast}\,f,\, \Lambda_{i}\,(\,x\,)\,P_{F_{\,i}\,(\,x\,)}\,V^{\,\ast}\,T\,\left(\,V^{\,-\, 1}\,\right)^{\,\ast}\,f\,\right>\,d\mu_{x}\\
&=\,\sum\limits_{i \,\in\, [\,m\,]}\,\int\limits_{\,X}\,v_{\,i}^{\,2}\,(\,x\,)\,\left<\,\Lambda_{i}\,(\,x\,)\,P_{F_{\,i}\,(\,x\,)}\,V^{\,\ast}\,\left(\,V^{\,-\, 1}\,\right)^{\,\ast}\,U\,f,\, \Lambda_{i}\,(\,x\,)\,P_{F_{\,i}\,(\,x\,)}\,V^{\,\ast}\,\left(\,V^{\,-\, 1}\,\right)^{\,\ast}\,T\,f\,\right>\,d\mu_{x}\\
&=\,\sum\limits_{i \,\in\, [\,m\,]}\,\int\limits_{\,X}\,v_{\,i}^{\,2}\,(\,x\,)\,\left<\,\Lambda_{i}\,(\,x\,)\,P_{F_{\,i}\,(\,x\,)}\,U\,f,\, \Lambda_{i}\,(\,x\,)\,P_{F_{\,i}\,(\,x\,)}\,T\,f\,\right>\,d\mu_{x}.
\end{align*}
On the other hand, for each \,$f \,\in\, H$, it is easy to verify that
\[\sum\limits_{i \,\in\, [\,m\,]}\,\int\limits_{\,X}\,v_{\,i}^{\,2}\,(\,x\,)\,\left<\,\Lambda_{i}\,(\,x\,)\,P_{F_{\,i}\,(\,x\,)}\,U\,f,\, \Lambda_{i}\,(\,x\,)\,P_{F_{\,i}\,(\,x\,)}\,T\,f\,\right>\,d\mu_{x}\leq\,B\, \left\|\,V^{\,-\, 1}\,\right\|^{\,2}\,\|\,f\,\|^{\,2}.\]
This completes the proof.
\end{proof}

Next, we will see that the intersection of components of a W. C. C. K. G. F. F. with a closed subspace is a W. C. C. K. G. F. F. for the smaller space.

\begin{theorem}
Let \,$\left\{\,F\,(\,x\,),\, \Lambda\,(\,x\,),\, v\,(\,x\,)\,\right\}_{x \,\in\, X}$\, and \,$\left\{\,G\,(\,x\,),\, \Gamma\,(\,x\,),\, w\,(\,x\,)\,\right\}_{x \,\in\, X}$\, be W. C. C. K. G. F. F. for \,$H$\, and \,$W$\, be a closed subspace of \,$H$.\,Then the families \,$\left\{\,F\,(\,x\,)\,\cap\,W,\, \Lambda\,(\,x\,),\, v\,(\,x\,)\,\right\}_{x \,\in\, X}$\, and \,$\left\{\,G\,(\,x\,)\,\cap\,W,\, \Gamma\,(\,x\,),\, w\,(\,x\,)\,\right\}_{x \,\in\, X}$\, are W. C. C. K. G. F. F. for \,$W$.  
\end{theorem}

\begin{proof}
The operators \,$P_{F\,(\,x\,)\,\cap\,W} \,=\, P_{F\,(\,x\,)}\left(\,P_{W}\,\right)$\, and \,$P_{G\,(\,x\,)\,\cap\,W} \,=\, P_{G\,(\,x\,)}\left(\,P_{W}\,\right)$\, are orthogonal projections of \,$H$\, onto \,$F\,(\,x\,)\,\cap\,W$\, and \,$G\,(\,x\,)\,\cap\,W$, respectively.\,Let \,$\sigma$\, be a measurable subset of \,$X$.\,Then for each \,$f \,\in\, H$, we have
\begin{align*}
&\int\limits_{\,\sigma}\,v^{\,2}\,(\,x\,)\,\left<\,\Lambda\,(\,x\,)\,P_{\,F\,(\,x\,)}\,U\,f,\, \Lambda\,(\,x\,)\,P_{\,F\,(\,x\,)}\,T\,f\,\right>\,d\mu_{x} \,+\, \\
&\hspace{1cm}\int\limits_{\,\sigma^{\,c}}\,w^{\,2}\,(\,x\,)\,\left<\,\Gamma\,(\,x\,)\,P_{\,G\,(\,x\,)}\,U\,f,\, \Gamma\,(\,x\,)\,P_{\,G\,(\,x\,)}\,T\,f\,\right>\,d\mu_{x}\\
&=\,\int\limits_{\,\sigma}\,v^{\,2}\,(\,x\,)\,\left<\,\Lambda\,(\,x\,)\,P_{\,F\,(\,x\,)}\,P_{W}\,U\,f,\, \Lambda\,(\,x\,)\,P_{\,F\,(\,x\,)}\,P_{W}\,T\,f\,\right>\,d\mu_{x} \,+\, \\
&\hspace{1cm}\int\limits_{\,\sigma^{\,c}}\,w^{\,2}\,(\,x\,)\,\left<\,\Gamma\,(\,x\,)\,P_{\,G\,(\,x\,)}\,P_{W}\,U\,f,\, \Gamma\,(\,x\,)\,P_{\,G\,(\,x\,)}\,P_{W}\,T\,f\,\right>\,d\mu_{x}\\
&=\,\int\limits_{\,\sigma}\,v^{\,2}\,(\,x\,)\,\left<\,\Lambda\,(\,x\,)\,P_{\,F\,(\,x\,)\,\cap\,W}\,U\,f,\, \Lambda\,(\,x\,)\,P_{\,F\,(\,x\,)\,\cap\,W}\,T\,f\,\right>\,d\mu_{x} \,+\, \\
&\hspace{1cm}\int\limits_{\,\sigma^{\,c}}\,w^{\,2}\,(\,x\,)\,\left<\,\Gamma\,(\,x\,)\,P_{\,G\,(\,x\,)\,\cap\,W}\,U\,f,\, \Gamma\,(\,x\,)\,P_{\,G\,(\,x\,)\,\cap\,W}\,T\,f\,\right>\,d\mu_{x}.
\end{align*} 
This completes the proof.
\end{proof}

The following theorem states the equivalence between W. C. C. K. G. F. F. and a bounded linear operator. 

\begin{theorem}
Let \,$K,\, V \,\in\, \mathcal{B}\,(\,H\,)$\, be an invertible operators such that \,$V^{\,\ast}$\, commutes with \,$T,\,U$.\,Let \,$\Lambda_{T\,U} \,=\, \left\{\,\left(\,F_{\,i}\,(\,x\,),\, \Lambda_{i}\,(\,x\,),\, v_{\,i}\,(\,x\,)\,\right)\,\right\}_{i \,\in\, [\,m\,],\, x \,\in\, \sigma_{\,i}}$\, be a W. C. C. K. G. F. F. for \,$H$\, with universal bounds \,$A$\, and \,$B$.\,Then the families \,$\Delta_{T\,U} \,=\, \left\{\,\left(\,V\,F_{\,i}\,(\,x\,),\, \Lambda_{i}\,(\,x\,)\,P_{F_{\,i}\,(\,x\,)}\,V^{\,\ast},\, v_{\,i}\,(\,x\,)\,\right)\,\right\}_{i \,\in\, [\,m\,],\, x \,\in\, \sigma_{\,i}}$\, is a W. C. C. K. G. F. F. for \,$H$\, if and only if there exists a \,$\delta \,>\, 0$\, such that for each \,$f \,\in\, H$, we have \,$\left\|\,V^{\,\ast}\,f\,\right\| \,\geq\, \delta\,\left\|\,K^{\,\ast}\,f\,\right\|$.   
\end{theorem}

\begin{proof}
Suppose that \,$\Delta_{T\,U}$\, is a W. C. C. K. G. F. F. for \,$H$\, with bounds \,$C$\, and \,$D$.\,Then for each \,$f \,\in\, H$, using the Theorem \ref{th1.01}, we have
\begin{align*}
&C\,\left\|\,K^{\,\ast}\,f\,\right\|^{\,2}\\
& \,\leq\, \sum\limits_{i \,\in\, [\,m\,]}\,\int\limits_{\,X}\,v_{\,i}^{\,2}\,(\,x\,)\,\left<\,\Lambda_{i}\,(\,x\,)\,P_{F_{\,i}\,(\,x\,)}\,V^{\,\ast}\,P_{V\,F_{\,i}\,(\,x\,)}\,U\,f,\, \Lambda_{i}\,(\,x\,)\,P_{F_{\,i}\,(\,x\,)}\,V^{\,\ast}\,P_{V\,F_{\,i}\,(\,x\,)}\,T\,f\,\right>\,d\mu_{x}\\
&=\,\sum\limits_{i \,\in\, [\,m\,]}\,\int\limits_{\,X}\,v_{\,i}^{\,2}\,(\,x\,)\,\left<\,\Lambda_{i}\,(\,x\,)\,P_{F_{\,i}\,(\,x\,)}\,V^{\,\ast}\,U\,f,\, \Lambda_{i}\,(\,x\,)\,P_{F_{\,i}\,(\,x\,)}\,V^{\,\ast}\,T\,f\,\right>\,d\mu_{x}\\
&=\,\sum\limits_{i \,\in\, [\,m\,]}\,\int\limits_{\,X}\,v_{\,i}^{\,2}\,(\,x\,)\,\left<\,\Lambda_{i}\,(\,x\,)\,P_{F_{\,i}\,(\,x\,)}\,U\,V^{\,\ast}\,f,\, \Lambda_{i}\,(\,x\,)\,P_{F_{\,i}\,(\,x\,)}\,T\,V^{\,\ast}\,f\,\right>\,d\mu_{x}\\
&\leq\,B\,\left\|\,V^{\,\ast}\,f\,\right\|^{\,2}. 
\end{align*}
Thus, 
\[\left\|\,V^{\,\ast}\,f\,\right\| \,\geq\, \sqrt{C \,/\, B}\,\left\|\,K^{\,\ast}\,f\,\right\|\; \;\forall\; f \,\in\, H.\]
Conversely, suppose \,$\left\|\,V^{\,\ast}\,f\,\right\| \,\geq\, \delta\,\left\|\,K^{\,\ast}\,f\,\right\|\; \;\forall\; f \,\in\, H$.\,Now, for \,$f \,\in\, H$, we have
\begin{align*}
&\sum\limits_{i \,\in\, [\,m\,]}\,\int\limits_{\,X}\,v_{\,i}^{\,2}\,(\,x\,)\,\left<\,\Lambda_{i}\,(\,x\,)\,P_{F_{\,i}\,(\,x\,)}\,V^{\,\ast}\,P_{V\,F_{\,i}\,(\,x\,)}\,U\,f,\, \Lambda_{i}\,(\,x\,)\,P_{F_{\,i}\,(\,x\,)}\,V^{\,\ast}\,P_{V\,F_{\,i}\,(\,x\,)}\,T\,f\,\right>\,d\mu_{x}\\
&=\,\sum\limits_{i \,\in\, [\,m\,]}\,\int\limits_{\,X}\,v_{\,i}^{\,2}\,(\,x\,)\,\left<\,\Lambda_{i}\,(\,x\,)\,P_{F_{\,i}\,(\,x\,)}\,U\,V^{\,\ast}\,f,\, \Lambda_{i}\,(\,x\,)\,P_{F_{\,i}\,(\,x\,)}\,T\,V^{\,\ast}\,f\,\right>\,d\mu_{x}\\
&\geq\, A\,\left\|\,K^{\,\ast}\,V^{\,\ast}\,f\,\right\|^{\,2} \,\geq\, A\,\left\|\,K^{\,-\, 1}\,\right\|^{\,-\, 2}\,\left\|\,V^{\,\ast}\,f\,\right\|^{\,2} \,\geq\, A\,\delta^{\,2}\,\left\|\,K^{\,-\, 1}\,\right\|^{\,-\, 2}\,\left\|\,K^{\,\ast}\,f\,\right\|^{\,2}.
\end{align*}
This completes the proof.   
\end{proof}

The next theorem shows that it is enough to cheek continuous weaving controlled \,$K$-$g$-fusion woven on smaller measurable space than the original.

\begin{theorem}
Let \,$\left\{\,\left(\,F_{\,i}\,(\,x\,),\, \Lambda_{i}\,(\,x\,),\, v_{\,i}\,(\,x\,)\,\right)\,\right\}_{x \,\in\, X}$\, be a continuous \,$(\,T,\,U\,)$-controlled $K$-$g$-fusion frame for \,$H$\, with universal bounds \,$A_{\,i}$\, and \,$B_{\,i}$, for each \,$i \,\in\, [\,m\,]$.\,If there exists a measurable subset \,$Y \,\subset\, X$\, such that the family of continuous \,$(\,T,\,U\,)$-controlled $K$-$g$-fusion frame \,$\left\{\,\left(\,F_{\,i}\,(\,x\,),\, \Lambda_{i}\,(\,x\,),\, v_{\,i}\,(\,x\,)\,\right)\,\right\}_{i \,\in\, [\,m\,],\, x \,\in\, Y}$\, is a W. C. C. K. G. F. F. for \,$H$\, with universal frame bounds \,$A$\, and \,$B$.\,Then \,$\left\{\,\left(\,F_{\,i}\,(\,x\,),\, \Lambda_{i}\,(\,x\,),\, v_{\,i}\,(\,x\,)\,\right)\,\right\}_{i \,\in\, [\,m\,],\, x \,\in\, X}$\, is a W. C. C. K. G. F. F. for \,$H$\, with universal frame bounds \,$A$\, and \,$\sum\limits_{i \,\in\, [\,m\,]}\,B_{\,i}$.      
\end{theorem} 

\begin{proof}
Let \,$\left\{\,\sigma_{\,i}\,\right\}_{i \,\in\, [\,m\,]}$\, be an arbitrary partition of \,$X$.\,For each \,$f \,\in\, H$, we define \,$\varphi \,:\, X \,\to\, \mathbb{C}$\, by 
\[\varphi\,(\,x\,) \,=\, \sum\limits_{i \,\in\, [\,m\,]}\,\chi_{\sigma_{\,i}}\,(\,x\,)\,\left<\,\Lambda_{i}\,(\,x\,)\,P_{F_{\,i}\,(\,x\,)}\,U\,f,\, \Lambda_{i}\,(\,x\,)\,P_{F_{\,i}\,(\,x\,)}\,T\,f\,\right>.\] 
Then \,$\varphi$\, is measurable.\,Now, for each \,$f \,\in\, H$, we have
\begin{align*}
&\sum\limits_{i \,\in\, [\,m\,]}\,\int\limits_{\,\sigma_{\,i}}\,v_{\,i}^{\,2}\,(\,x\,)\,\left<\,\Lambda_{i}\,(\,x\,)\,P_{F_{\,i}\,(\,x\,)}\,U\,f,\, \Lambda_{i}\,(\,x\,)\,P_{F_{\,i}\,(\,x\,)}\,T\,f\,\right>\,d\mu_{x}\\
&\leq\,\sum\limits_{i \,\in\, [\,m\,]}\,\int\limits_{\,X}\,v_{\,i}^{\,2}\,(\,x\,)\,\left<\,\Lambda_{i}\,(\,x\,)\,P_{F_{\,i}\,(\,x\,)}\,U\,f,\, \Lambda_{i}\,(\,x\,)\,P_{F_{\,i}\,(\,x\,)}\,T\,f\,\right>\,d\mu_{x}\\
&\leq\,\left(\,\sum\limits_{i \,\in\, [\,m\,]}\,B_{\,i}\,\right)\,\|\,f\,\|^{\,2}.
\end{align*}
It is easy to verify that \,$\left\{\,\sigma_{\,i}\,\cap\,Y\,\right\}_{i \,\in\, [\,m\,]}$\, is a partitions of \,$Y$.\,Thus, the family \,$\left\{\,\left(\,F_{\,i}\,(\,x\,),\, \Lambda_{i}\,(\,x\,),\, v_{\,i}\,(\,x\,)\,\right)\,\right\}_{i \,\in\, [\,m\,],\, x \,\in\, \sigma_{\,i}\,\cap\,Y}$\, is a continuous \,$(\,T,\,U\,)$-controlled $K$-$g$-fusion frame for \,$H$\, with lowest frame bounds \,$A$.\,Therefore 
\begin{align*}
&\sum\limits_{i \,\in\, [\,m\,]}\,\int\limits_{\,\sigma_{\,i}}\,v_{\,i}^{\,2}\,(\,x\,)\,\left<\,\Lambda_{i}\,(\,x\,)\,P_{F_{\,i}\,(\,x\,)}\,U\,f,\, \Lambda_{i}\,(\,x\,)\,P_{F_{\,i}\,(\,x\,)}\,T\,f\,\right>\,d\mu_{x}\\
&\geq\,\sum\limits_{i \,\in\, [\,m\,]}\,\int\limits_{\,\sigma_{\,i}\,\cap\,Y}\,v_{\,i}^{\,2}\,(\,x\,)\,\left<\,\Lambda_{i}\,(\,x\,)\,P_{F_{\,i}\,(\,x\,)}\,U\,f,\, \Lambda_{i}\,(\,x\,)\,P_{F_{\,i}\,(\,x\,)}\,T\,f\,\right>\,d\mu_{x}\\
&\geq\,A\,\left\|\,K^{\,\ast}\,f\,\right\|^{\,2}.
\end{align*} 
This completes the proof.
\end{proof}

In the following theorem, we shows that it is possible to remove vectors from continuous controlled $K$-$g$-fusion frames and still be left with woven frames.

\begin{theorem}
Let \,$\left\{\,\left(\,F_{\,i}\,(\,x\,),\, \Lambda_{i}\,(\,x\,),\, v_{\,i}\,(\,x\,)\,\right)\,\right\}_{i \,\in\, [\,m\,],\, x \,\in\, \sigma_{\,i}}$\, be a W. C. C. K. G. F. F. for \,$H$\, with universal bounds \,$A$\, and \,$B$.\,If there exists \,$0 \,<\, D \,<\, A$\, and a measurable subset \,$Y \,\subset\, X$\, and \,$n \,\in\, [\,m\,]$\, such that for \,$f \,\in\, H$ 
\[\sum\limits_{i \,\in\, [\,m\,] \,\setminus\, \{\,n\,\}}\,\int\limits_{\,X \,\setminus\, Y}\,v_{\,i}^{\,2}\,(\,x\,)\,\left<\,\Lambda_{i}\,(\,x\,)\,P_{F_{\,i}\,(\,x\,)}\,U\,f,\, \Lambda_{i}\,(\,x\,)\,P_{F_{\,i}\,(\,x\,)}\,T\,f\,\right>\,d\mu_{x} \,\leq\, D\,\left\|\,K^{\,\ast}\,f\,\right\|^{\,2}.\]
Then the family \,$\left\{\,\left(\,F_{\,i}\,(\,x\,),\, \Lambda_{i}\,(\,x\,),\, v_{\,i}\,(\,x\,)\,\right)\,\right\}_{i \,\in\, [\,m\,],\, x \,\in\, Y}$\, is a W. C. C. K. G. F. F. for \,$H$\, with frame bounds \,$A \,-\, D$\, and \,$B$.
\end{theorem} 

\begin{proof}
Suppose that \,$\left\{\,\sigma_{\,i}\,\right\}_{i \,\in\, [\,m\,]}$\, and \,$\left\{\,\gamma_{\,i}\,\right\}_{i \,\in\, [\,m\,]}$\, are partitions of \,$Y$\, and \,$X \,\setminus\, Y$, respectively.\,For a given \,$f \,\in\, H$, we define \,$\varphi \,:\, Y \,\to\, \mathbb{C}$\, by 
\[\varphi\,(\,x\,) \,=\, \sum\limits_{i \,\in\, [\,m\,]}\,\chi_{\sigma_{\,i}}\,(\,x\,)\,\left<\,\Lambda_{i}\,(\,x\,)\,P_{F_{\,i}\,(\,x\,)}\,U\,f,\, \Lambda_{i}\,(\,x\,)\,P_{F_{\,i}\,(\,x\,)}\,T\,f\,\right>\]
and \,$\phi \,:\, X \,\to\, \mathbb{C}$\, by 
\[\phi\,(\,x\,) \,=\, \sum\limits_{i \,\in\, [\,m\,]}\,\chi_{\sigma_{\,i}\,\cup\,\gamma_{\,i}}\,(\,x\,)\,\left<\,\Lambda_{i}\,(\,x\,)\,P_{F_{\,i}\,(\,x\,)}\,U\,f,\, \Lambda_{i}\,(\,x\,)\,P_{F_{\,i}\,(\,x\,)}\,T\,f\,\right>.\]
Since \,$\left\{\,\left(\,F_{\,i}\,(\,x\,),\, \Lambda_{i}\,(\,x\,),\, v_{\,i}\,(\,x\,)\,\right)\,\right\}_{i \,\in\, [\,m\,],\, x \,\in\, \sigma_{\,i}\,\cup\,\gamma_{\,i}}$\, is a continuous \,$(\,T,\,U\,)$-controlled $K$-$g$-fusion frame for \,$H$\, and \,$\varphi \,=\, \phi\,|_{\,Y}$, \,$\varphi$\, and \,$\phi$\, are measurable.\,So, for each \,$f \,\in\, H$, we have
\begin{align*}
&\sum\limits_{i \,\in\, [\,m\,]}\,\int\limits_{\,\sigma_{\,i}}\,v_{\,i}^{\,2}\,(\,x\,)\,\left<\,\Lambda_{i}\,(\,x\,)\,P_{F_{\,i}\,(\,x\,)}\,U\,f,\, \Lambda_{i}\,(\,x\,)\,P_{F_{\,i}\,(\,x\,)}\,T\,f\,\right>\,d\mu_{x}\\
&\leq\, \sum\limits_{i \,\in\, [\,m\,]}\,\int\limits_{\,\sigma_{\,i}\,\cup\,\gamma_{\,i}}\,v_{\,i}^{\,2}\,(\,x\,)\,\left<\,\Lambda_{i}\,(\,x\,)\,P_{F_{\,i}\,(\,x\,)}\,U\,f,\, \Lambda_{i}\,(\,x\,)\,P_{F_{\,i}\,(\,x\,)}\,T\,f\,\right>\,d\mu_{x} \,\leq\, B\,\|\,f\,\|^{\,2}.
\end{align*} 
Now, we assume that \,$\left\{\,\xi_{\,i}\,\right\}_{i \,\in\, [\,m\,]}$\, such that \,$\xi_{\,n} \,=\, \theta$.\,Then \,$\left\{\,\xi_{\,i} \,\cup\, \sigma_{\,i}\,\right\}_{i \,\in\, [\,m\,]}$\, is a partition of \,$X$\, and so for any \,$f \,\in\, H$, we have
\begin{align*}
&\sum\limits_{i \,\in\, [\,m\,]}\,\int\limits_{\,\sigma_{\,i}}\,v_{\,i}^{\,2}\,(\,x\,)\,\left<\,\Lambda_{i}\,(\,x\,)\,P_{F_{\,i}\,(\,x\,)}\,U\,f,\, \Lambda_{i}\,(\,x\,)\,P_{F_{\,i}\,(\,x\,)}\,T\,f\,\right>\,d\mu_{x}\\
&=\,\sum\limits_{i \,\in\, [\,m\,] \,\setminus\, \{\,n\,\}}\,\bigg[\int\limits_{\,\xi_{\,i}\,\cup\,\sigma_{\,i}}\,v_{\,i}^{\,2}\,(\,x\,)\,\left<\,\Lambda_{i}\,(\,x\,)\,P_{F_{\,i}\,(\,x\,)}\,U\,f,\, \Lambda_{i}\,(\,x\,)\,P_{F_{\,i}\,(\,x\,)}\,T\,f\,\right>\,d\mu_{x}\,-\\
&\hspace{1cm}\int\limits_{\,\xi_{\,i}}\,v_{\,i}^{\,2}\,(\,x\,)\,\left<\,\Lambda_{i}\,(\,x\,)\,P_{F_{\,i}\,(\,x\,)}\,U\,f,\, \Lambda_{i}\,(\,x\,)\,P_{F_{\,i}\,(\,x\,)}\,T\,f\,\right>\,d\mu_{x}\,+\\
&\hspace{1cm}\int\limits_{\,\sigma_{\,n}}\,v_{\,i}^{\,2}\,(\,x\,)\,\left<\,\Lambda_{i}\,(\,x\,)\,P_{F_{\,i}\,(\,x\,)}\,U\,f,\, \Lambda_{i}\,(\,x\,)\,P_{F_{\,i}\,(\,x\,)}\,T\,f\,\right>\,d\mu_{x}\,\bigg]\\
&\geq\, \sum\limits_{i \,\in\, [\,m\,] \,\setminus\, \{\,n\,\}}\,\bigg[\int\limits_{\,\xi_{\,i}\,\cup\,\sigma_{\,i}}\,v_{\,i}^{\,2}\,(\,x\,)\,\left<\,\Lambda_{i}\,(\,x\,)\,P_{F_{\,i}\,(\,x\,)}\,U\,f,\, \Lambda_{i}\,(\,x\,)\,P_{F_{\,i}\,(\,x\,)}\,T\,f\,\right>\,d\mu_{x}\,-\\
&\hspace{1cm}\int\limits_{\,X \,\setminus\, Y}\,v_{\,i}^{\,2}\,(\,x\,)\,\left<\,\Lambda_{i}\,(\,x\,)\,P_{F_{\,i}\,(\,x\,)}\,U\,f,\, \Lambda_{i}\,(\,x\,)\,P_{F_{\,i}\,(\,x\,)}\,T\,f\,\right>\,d\mu_{x}\,+\\
&\hspace{1cm}\int\limits_{\,\sigma_{\,n}}\,v_{\,i}^{\,2}\,(\,x\,)\,\left<\,\Lambda_{i}\,(\,x\,)\,P_{F_{\,i}\,(\,x\,)}\,U\,f,\, \Lambda_{i}\,(\,x\,)\,P_{F_{\,i}\,(\,x\,)}\,T\,f\,\right>\,d\mu_{x}\,\bigg]\\
&=\,\sum\limits_{i \,\in\, [\,m\,]}\,\int\limits_{\,\xi_{\,i}\,\cup\,\sigma_{\,i}}\,v_{\,i}^{\,2}\,(\,x\,)\,\left<\,\Lambda_{i}\,(\,x\,)\,P_{F_{\,i}\,(\,x\,)}\,U\,f,\, \Lambda_{i}\,(\,x\,)\,P_{F_{\,i}\,(\,x\,)}\,T\,f\,\right>\,d\mu_{x} \,-\\
&\hspace{1cm}\sum\limits_{i \,\in\, [\,m\,] \,\setminus\, \{\,n\,\}}\,\int\limits_{\,X \,\setminus\, Y}\,v_{\,i}^{\,2}\,(\,x\,)\,\left<\,\Lambda_{i}\,(\,x\,)\,P_{F_{\,i}\,(\,x\,)}\,U\,f,\, \Lambda_{i}\,(\,x\,)\,P_{F_{\,i}\,(\,x\,)}\,T\,f\,\right>\,d\mu_{x}\\
&\geq\, (\,A \,-\, D\,)\,\left\|\,K^{\,\ast}\,f\,\right\|^{\,2}.
\end{align*} 
This completes the proof.    
\end{proof}
 
\begin{proposition} 
Let \,$K \,\in\, \mathcal{B}\,(\,H\,)$\, be an invertible operator, \,$V \,\in\, \mathcal{B}\,(\,H\,)$\, be a unitary operator and \,$\left\{\,\left(\,F\,(\,x\,),\, \Lambda\,(\,x\,),\, v\,(\,x\,)\,\right)\,\right\}_{x \,\in\, X}$\, be a continuous \,$(\,T,\,U\,)$-controlled \,$K$-$g$-fusion frame for \,$H$\, with bounds \,$A,\,B$.\,If \,$\left\|\,I_{H} \,-\, V\,\right\|^{\,2}\,\left\|\,K^{\,-\, 1}\,\right\|^{\,2} \,\leq\, A \,/\, B$\, and \,$V$\, commutes with \,$T,\,U$, then \,$\Lambda \,=\, \left\{\,\left(\,F\,(\,x\,),\, \Lambda\,(\,x\,),\, v\,(\,x\,)\,\right)\,\right\}_{x \,\in\, X}$\, and \,$\Lambda^{\,\prime} \,=\, \left\{\,\left(\,V^{\,-\, 1}\,F\,(\,x\,),\, \Lambda\,(\,x\,)\,V,\, v\,(\,x\,)\,\right)\,\right\}_{x \,\in\, X}$\, are W. C. C. K. G. F. F. for \,$H$.    
\end{proposition}
 
\begin{proof}
Let \,$\sigma$\, be a partition of \,$X$.\,Since \,$K \,\in\, \mathcal{B}\,(\,H\,)$\, is an invertible operator, for \,$f \,\in\, H$, we have \,$\|\,f\,\|^{\,2} \,\leq\, \left\|\,K^{\,-\, 1}\,\right\|^{\,2}\,\left\|\,K^{\,\ast}\,f\,\right\|^{\,2}$.\,Now, for each \,$f \,\in\, H$, we have
\begin{align*}
&\int\limits_{\,\sigma}\,v^{\,2}\,(\,x\,)\,\left<\,\Lambda\,(\,x\,)\,P_{\,F\,(\,x\,)}\,U\,f,\, \Lambda\,(\,x\,)\,P_{\,F\,(\,x\,)}\,T\,f\,\right>\,d\mu_{x} \,+\, \\
&\hspace{1cm}\int\limits_{\,\sigma^{\,c}}\,v^{\,2}\,(\,x\,)\,\left<\,\Lambda\,(\,x\,)\,V\,P_{\,V^{\,-\, 1}\,F\,(\,x\,)}\,U\,f,\, \Lambda\,(\,x\,)\,V\,P_{\,V^{\,-\, 1}\,F\,(\,x\,)}\,T\,f\,\right>\,d\mu_{x}\\
&=\,\int\limits_{\,\sigma}\,v^{\,2}\,(\,x\,)\,\left<\,\Lambda\,(\,x\,)\,P_{\,F\,(\,x\,)}\,U\,f,\, \Lambda\,(\,x\,)\,P_{\,F\,(\,x\,)}\,T\,f\,\right>\,d\mu_{x} \,+\, \\
&\hspace{1cm}\int\limits_{\,\sigma^{\,c}}\,v^{\,2}\,(\,x\,)\,\left<\,\Lambda\,(\,x\,)\,P_{\,F\,(\,x\,)}\,U\,V\,f,\, \Lambda\,(\,x\,)\,P_{\,F\,(\,x\,)}\,T\,V\,f\,\right>\,d\mu_{x}\\
&\geq\,\int\limits_{\,X}\,v^{\,2}\,(\,x\,)\,\left<\,\Lambda\,(\,x\,)\,P_{\,F\,(\,x\,)}\,U\,f,\, \Lambda\,(\,x\,)\,P_{\,F\,(\,x\,)}\,T\,f\,\right>\,d\mu_{x} \,-\, \\
&\int\limits_{\,\sigma^{\,c}}\,v^{\,2}\,(\,x\,)\,\left<\,\Lambda\,(\,x\,)\,P_{\,F\,(\,x\,)}\,U\,\left(\,I_{H} \,-\, V\,\right)\,f,\, \Lambda\,(\,x\,)\,P_{\,F\,(\,x\,)}\,T\,\left(\,I_{H} \,-\, V\,\right)\,f\,\right>\,d\mu_{x}\\
&\geq\,A\,\left\|\,K^{\,\ast}\,f\,\right\|^{\,2} \,-\, B\,\left\|\,I_{H} \,-\, V\,\right\|^{\,2}\,\|\,f\,\|^{\,2}\\
&\geq\,A\,\left\|\,K^{\,\ast}\,f\,\right\|^{\,2} \,-\, B\,\left\|\,I_{H} \,-\, V\,\right\|^{\,2}\,\left\|\,K^{\,-\, 1}\,\right\|^{\,2}\,\left\|\,K^{\,\ast}\,f\,\right\|^{\,2}\\
&=\,\left(\,A \,-\, B\,\left\|\,I_{H} \,-\, V\,\right\|^{\,2}\,\left\|\,K^{\,-\, 1}\,\right\|^{\,2}\,\right)\,\left\|\,K^{\,\ast}\,f\,\right\|^{\,2}. 
\end{align*}
Hence, the families \,$\Lambda$\, and \,$\Lambda^{\,\prime}$\, are W. C. C. K. G. F. F. for \,$H$.
\end{proof} 

Next, we will see that under some sufficient conditions sum of two continuous \,$(\,T,\,U\,)$-controlled \,$K$-$g$-fusion frames is woven with itself. 

\begin{theorem}
Let \,$K \,\in\, \mathcal{B}\,(\,H\,)$\, be an invertible operator, the families \,$\Lambda \,=\, \left\{\,\left(\,F\,(\,x\,),\, \Lambda\,(\,x\,),\, v\,(\,x\,)\,\right)\,\right\}_{x \,\in\, X}$\,  and \,$\Gamma \,=\, \left\{\,\left(\,G\,(\,x\,),\, \Lambda\,(\,x\,),\, v\,(\,x\,)\,\right)\,\right\}_{x \,\in\, X}$\, be continuous \,$(\,T,\,U\,)$-controlled \,$K$-$g$-fusion frame for \,$H$\, with bounds \,$A,\,B$\, and \,$C,\,D$, respectively.\,Suppose for each \,$x \,\in\, X$,
\begin{itemize}
\item[$(i)$]$F\,(\,x\,) \,\subset\, G\,(\,x\,)^{\,\perp}$,
\item[$(ii)$]$\Lambda\,(\,x\,)\,P_{\,F\,(\,x\,)}\,\mathcal{R}\,(\,U\,) \,\perp\, \Lambda\,(\,x\,)\,P_{\,G\,(\,x\,)}\,\mathcal{R}\,(\,T\,)$,
\item[$(iii)$]$\Lambda\,(\,x\,)\,P_{\,F\,(\,x\,)}\,\mathcal{R}\,(\,T\,) \,\perp\, \Lambda\,(\,x\,)\,P_{\,G\,(\,x\,)}\,\mathcal{R}\,(\,U\,)$.
\end{itemize}
If for any partition \,$\sigma$\, of \,$X$, \,$\left(\,T_{\Gamma}^{\,\sigma}\,\right)^{\,\ast}$\, is bounded below then the families \,$\Delta \,=\, \left\{\,\left(\,F\,(\,x\,) \,+\, G\,(\,x\,),\, \Lambda\,(\,x\,),\, v\,(\,x\,)\,\right)\,\right\}_{x \,\in\, X}$\, and \,$\Lambda$\, are W. C. C. K. G. F. F. for \,$H$.  
\end{theorem}

\begin{proof}
Since for each \,$x \,\in\, X$, \,$F\,(\,x\,) \,\subset\, G\,(\,x\,)^{\,\perp}$, we have \,$P_{\,F\,(\,x\,) \,+\, G\,(\,x\,)} \,=\, P_{\,F\,(\,x\,)} \,+\, P_{\,F\,(\,x\,)}$.\,Now, for each \,$x \,\in\, X$, using the given conditions $(ii)$\, and $(iii)$, we have
\begin{align}
&\int\limits_{\,X}\,v^{\,2}\,(\,x\,)\,\left<\,\Lambda\,(\,x\,)\,P_{\,F\,(\,x\,) \,+\, G\,(\,x\,)}\,U\,f,\, \Lambda\,(\,x\,)\,P_{\,F\,(\,x\,) \,+\, G\,(\,x\,)}\,T\,f\,\right>\,d\mu_{x}\nonumber\\
&=\,\int\limits_{\,X}\,v^{\,2}\,(\,x\,)\,\left<\,\Lambda\,(\,x\,)\,\left(\,P_{\,F\,(\,x\,)} \,+\, P_{\,G\,(\,x\,)}\,\right)\,U\,f,\, \Lambda\,(\,x\,)\,\left(\,P_{\,F\,(\,x\,)} \,+\, P_{\,G\,(\,x\,)}\,\right)\,T\,f\,\right>\,d\mu_{x}\nonumber\\
&=\,\int\limits_{\,X}\,v^{\,2}\,(\,x\,)\,\left<\,\Lambda\,(\,x\,)\,P_{\,F\,(\,x\,)}\,U\,f,\, \Lambda\,(\,x\,)\,P_{\,F\,(\,x\,)}\,T\,f\,\right>\,d\mu_{x} \,+\, \nonumber\\
&\hspace{1cm}\int\limits_{\,X}\,v^{\,2}\,(\,x\,)\,\left<\,\Lambda\,(\,x\,)\,P_{\,G\,(\,x\,)}\,U\,f,\, \Lambda\,(\,x\,)\,P_{\,G\,(\,x\,)}\,T\,f\,\right>\,d\mu_{x}\label{et1}\\
&\leq\,(\,B \,+\, D\,)\,\|\,f\,\|^{\,2}.\nonumber
\end{align}
On the other hand, from (\ref{et1}), we get
\begin{align*}
&\int\limits_{\,X}\,v^{\,2}\,(\,x\,)\,\left<\,\Lambda\,(\,x\,)\,P_{\,F\,(\,x\,) \,+\, G\,(\,x\,)}\,U\,f,\, \Lambda\,(\,x\,)\,P_{\,F\,(\,x\,) \,+\, G\,(\,x\,)}\,T\,f\,\right>\,d\mu_{x}\\
&\geq\,(\,A \,+\, C\,)\,\left\|\,K^{\,\ast}\,f\,\right\|^{\,2}\; \;\forall\; f \,\in\, H. 
\end{align*}
Thus, \,$\Delta$\, is a continuous \,$(\,T,\,U\,)$-controlled \,$K$-$g$-fusion frame for \,$H$\, with bound \,$(\,A \,+\, C\,)$\, and \,$(\,B \,+\, D\,)$.\\

Furthermore, since \,$K$\, is a invertible operator and for any partition \,$\sigma$\, of \,$X$, \,$\left(\,T_{\Gamma}^{\,\sigma}\,\right)^{\,\ast}$\, is bounded below, for each \,$f \,\in\, H$, there exists \,$M \,>\, 0$\, such that
\[\left\|\,\left(\,T_{\Gamma}^{\,\sigma}\,\right)^{\,\ast}\,f\,\right\|^{\,2} \,\geq\, M^{\,2}\,\|\,f\,\|^{\,2} \,\geq\, \dfrac{M^{\,2}}{\left\|\,K\,\right\|^{\,2}}\,\left\|\,K^{\,\ast}\,f\,\right\|^{\,2}.\]
Now, for each \,$f \,\in\, H$, we have
\begin{align*}
&\int\limits_{\,\sigma}\,v^{\,2}\,(\,x\,)\,\left<\,\Lambda\,(\,x\,)\,P_{\,F\,(\,x\,) \,+\, G\,(\,x\,)}\,U\,f,\, \Lambda\,(\,x\,)\,P_{\,F\,(\,x\,) \,+\, G\,(\,x\,)}\,T\,f\,\right>\,d\mu_{x} \,+\\
&\hspace{1cm}\int\limits_{\,\sigma^{\,c}}\,v^{\,2}\,(\,x\,)\,\left<\,\Lambda\,(\,x\,)\,P_{\,F\,(\,x\,)}\,U\,f,\, \Lambda\,(\,x\,)\,P_{\,F\,(\,x\,)}\,T\,f\,\right>\,d\mu_{x}\\
&=\,\int\limits_{\,X}\,v^{\,2}\,(\,x\,)\,\left<\,\Lambda\,(\,x\,)\,P_{\,F\,(\,x\,)}\,U\,f,\, \Lambda\,(\,x\,)\,P_{\,F\,(\,x\,)}\,T\,f\,\right>\,d\mu_{x} \,-\\
&\hspace{1cm}\,\int\limits_{\,\sigma}\,v^{\,2}\,(\,x\,)\,\left<\,\Lambda\,(\,x\,)\,P_{\,F\,(\,x\,)}\,U\,f,\, \Lambda\,(\,x\,)\,P_{\,F\,(\,x\,)}\,T\,f\,\right>\,d\mu_{x} \,+\\
&\hspace{.5cm}\,\int\limits_{\,\sigma}\,v^{\,2}\,(\,x\,)\,\left<\,\Lambda\,(\,x\,)\,\left(\,P_{\,F\,(\,x\,)} \,+\, P_{\,G\,(\,x\,)}\,\right)\,U\,f,\, \Lambda\,(\,x\,)\,\left(\,P_{\,F\,(\,x\,)} \,+\, P_{\,G\,(\,x\,)}\,\right)\,T\,f\,\right>\,d\mu_{x}\\
&=\,\int\limits_{\,X}\,v^{\,2}\,(\,x\,)\,\left<\,\Lambda\,(\,x\,)\,P_{\,F\,(\,x\,)}\,U\,f,\, \Lambda\,(\,x\,)\,P_{\,F\,(\,x\,)}\,T\,f\,\right>\,d\mu_{x} \,+\\
&\hspace{1cm}\,\int\limits_{\,\sigma}\,v^{\,2}\,(\,x\,)\,\left<\,\Lambda\,(\,x\,)\,P_{\,G\,(\,x\,)}\,U\,f,\, \Lambda\,(\,x\,)\,P_{\,G\,(\,x\,)}\,T\,f\,\right>\,d\mu_{x}\\
&\geq\, A\,\left\|\,K^{\,\ast}\,f\,\right\|^{\,2} \,+\, \left\|\,\left(\,T_{\Gamma}^{\,\sigma}\,\right)^{\,\ast}\,f\,\right\|^{\,2} \,\geq\, \left(\,A \,+\, \dfrac{M^{\,2}}{\left\|\,K\,\right\|^{\,2}}\,\right)\,\left\|\,K^{\,\ast}\,f\,\right\|^{\,2}. 
\end{align*}
On the other hand,
\begin{align*}
&\int\limits_{\,\sigma}\,v^{\,2}\,(\,x\,)\,\left<\,\Lambda\,(\,x\,)\,P_{\,F\,(\,x\,) \,+\, G\,(\,x\,)}\,U\,f,\, \Lambda\,(\,x\,)\,P_{\,F\,(\,x\,) \,+\, G\,(\,x\,)}\,T\,f\,\right>\,d\mu_{x} \,+\\
&\hspace{1cm}\int\limits_{\,\sigma^{\,c}}\,v^{\,2}\,(\,x\,)\,\left<\,\Lambda\,(\,x\,)\,P_{\,F\,(\,x\,)}\,U\,f,\, \Lambda\,(\,x\,)\,P_{\,F\,(\,x\,)}\,T\,f\,\right>\,d\mu_{x}\\
&\leq\,\int\limits_{\,X}\,v^{\,2}\,(\,x\,)\,\left<\,\Lambda\,(\,x\,)\,P_{\,F\,(\,x\,) \,+\, G\,(\,x\,)}\,U\,f,\, \Lambda\,(\,x\,)\,P_{\,F\,(\,x\,) \,+\, G\,(\,x\,)}\,T\,f\,\right>\,d\mu_{x} \,+\\
&\hspace{1cm}\int\limits_{\,X}\,v^{\,2}\,(\,x\,)\,\left<\,\Lambda\,(\,x\,)\,P_{\,F\,(\,x\,)}\,U\,f,\, \Lambda\,(\,x\,)\,P_{\,F\,(\,x\,)}\,T\,f\,\right>\,d\mu_{x}\\
&\leq\,(\,2\,B \,+\, D\,)\,\|\,f\,\|^{\,2}. 
\end{align*}
Thus, \,$\Delta$\, and \,$\Lambda$\, are W. C. C. K. G. F. F. for \,$H$.\,Similarly, it can be shown that \,$\Delta$\, and \,$\Gamma$\, are W. C. C. K. G. F. F. for \,$H$.\,This completes the proof.     
\end{proof}

In the following theorem, we present a sufficient condition for weaving continuous controlled $K$-$g$-fusion frame in terms of positive operators associated with given continuous controlled $K$-$g$-fusion frame.

\begin{theorem}
Let the families \,$\Lambda \,=\, \left\{\,\left(\,F\,(\,x\,),\, \Lambda\,(\,x\,),\, v\,(\,x\,)\,\right)\,\right\}_{x \,\in\, X}$\, and \,$\Gamma \,=\, \left\{\,\left(\,G\,(\,x\,),\, \Lambda\,(\,x\,),\, v\,(\,x\,)\,\right)\,\right\}_{x \,\in\, X}$\, be  continuous \,$(\,T,\,U\,)$-controlled \,$K$-$g$-fusion frame for \,$H$.\,Suppose for each \,$x \,\in\, X$, the operator \,$U_{\,x} \,:\, H \,\to\, H$\, defined by
\[\left<\,U_{\,x}\,(\,f\,),\, g\,\right>\,=\, \int\limits_{\,X}\,v^{\,2}\,(\,x\,)\,\left<\,T^{\,\ast}\,\Delta\,(\,x\,)\,U\,f,\, g\,\right>\,d\mu_{x},\; \;f,\, g \,\in\, H,\]
where \,$\Delta\,(\,x\,) \,=\, P_{\,G\,(\,x\,)}\,\Gamma^{\,\ast}\,(\,x\,)\,\Gamma\,(\,x\,)\,P_{\,G\,(\,x\,)} \,-\, P_{\,F\,(\,x\,)}\,\Lambda^{\,\ast}\,(\,x\,)\,\Lambda\,(\,x\,)\,P_{\,F\,(\,x\,)}$, is a positive operator.\,Then \,$\Lambda$\, and \,$\Gamma$\, are W. C. C. K. G. F. F. for \,$H$. 
\end{theorem}

\begin{proof}
Let \,$A,\, B$\, and \,$C,\, D$\, be frame bounds of \,$\Lambda$\, and \,$\Gamma$, respectively.\,Take \,$\sigma$\, be any partition of \,$X$.\,Then for each \,$f \,\in\, H$, we have
\begin{align*}
&A\,\left\|\,K^{\,\ast}\,f\,\right\|^{\,2} \,\leq\, \int\limits_{\,X}\,v^{\,2}\,(\,x\,)\,\left<\,\Lambda\,(\,x\,)\,P_{\,F\,(\,x\,)}\,U\,f,\, \Lambda\,(\,x\,)\,P_{\,F\,(\,x\,)}\,T\,f\,\right>\,d\mu_{x}\\
&=\,\int\limits_{\,\sigma}\,v^{\,2}\,(\,x\,)\,\left<\,\Lambda\,(\,x\,)\,P_{\,F\,(\,x\,)}\,U\,f,\, \Lambda\,(\,x\,)\,P_{\,F\,(\,x\,)}\,T\,f\,\right>\,d\mu_{x} \,+\,\\
&\hspace{1cm}\int\limits_{\,\sigma^{\,c}}\,v^{\,2}\,(\,x\,)\,\left<\,T^{\,\ast}\,P_{\,F\,(\,x\,)}\,\Lambda\,(\,x\,)^{\,\ast}\,\Lambda\,(\,x\,)\,P_{\,F\,(\,x\,)}\,U\,f,\, f\,\right>\,d\mu_{x}\\
&=\,\int\limits_{\,\sigma}\,v^{\,2}\,(\,x\,)\,\left<\,\Lambda\,(\,x\,)\,P_{\,F\,(\,x\,)}\,U\,f,\, \Lambda\,(\,x\,)\,P_{\,F\,(\,x\,)}\,T\,f\,\right>\,d\mu_{x} \,-\,\\
&\hspace{1cm}\int\limits_{\,\sigma^{\,c}}\,v^{\,2}\,(\,x\,)\,\left<\,T^{\,\ast}\,\Delta\,(\,x\,)\,U\,f,\, g\,\right>\,d\mu_{x} \,+\,\\
&\hspace{1cm} \int\limits_{\,\sigma^{\,c}}\,v^{\,2}\,(\,x\,)\,\left<\,T^{\,\ast}\,P_{\,G\,(\,x\,)}\,\Gamma\,(\,x\,)^{\,\ast}\,\Gamma\,(\,x\,)\,P_{\,G\,(\,x\,)}\,U\,f,\, f\,\right>\,d\mu_{x}\\
&\leq\, \int\limits_{\,\sigma}\,v^{\,2}\,(\,x\,)\,\left<\,\Lambda\,(\,x\,)\,P_{\,F\,(\,x\,)}\,U\,f,\, \Lambda\,(\,x\,)\,P_{\,F\,(\,x\,)}\,T\,f\,\right>\,d\mu_{x} \,+\,\\
&\hspace{1cm}\int\limits_{\,\sigma^{\,c}}\,v^{\,2}\,(\,x\,)\,\left<\,\Gamma\,(\,x\,)\,P_{\,G\,(\,x\,)}\,U\,f,\, \Gamma\,(\,x\,)\,P_{\,G\,(\,x\,)}\,T\,f\,\right>\,d\mu_{x}\\
&\leq\,(\,B \,+\, D\,)\,\|\,f\,\|^{\,2}.  
\end{align*}
Thus, \,$\Lambda$\, and \,$\Gamma$\, are W. C. C. K. G. F. F. for \,$H$\, with universal bounds \,$A$\, and \,$B \,+\, D$.   
\end{proof}

\begin{theorem}
Let \,$\left\{\,\left(\,F_{\,i}\,(\,x\,),\, \Lambda_{i}\,(\,x\,),\, v_{\,i}\,(\,x\,)\,\right)\,\right\}_{x \,\in\, X}$\, be a continuous \,$(\,T,\,U\,)$-controlled $K$-$g$-fusion frame for \,$H$\, with bounds \,$A_{i}$\, and \,$B_{i}$\, for each \,$i \,\in\, [\,m\,]$.\,Suppose \,$Y$\, be measurable subset \,$X$\, and there exists \,$N \,>\, 0$\, such that for all \,$i,\, k \,\in\, [\,m\,]$\, with \,$i \,\neq\, k$,
\[0 \,\leq\, \int\limits_{\,Y}\,\left<\,\Gamma_{i,\,k}\,U\,f,\, \Gamma_{i,\,k}\,T\,f\,\right>\,d\mu_{x} \,\leq\,N\,\min\,\{\,\Theta,\,\Omega\}, \,f \,\in\, H \]
where 
\[\Gamma_{i,\,k} \,=\, v_{\,i}^{\,2}\,(\,x\,)\,\Lambda_{i}\,(\,x\,)\,P_{F_{\,i}\,(\,x\,)}\,-\, v_{\,k}^{\,2}\,(\,x\,)\,\Lambda_{k}\,(\,x\,)\,P_{F_{\,k}\,(\,x\,)},\] 
\[\Theta \,=\, \int\limits_{\,Y}\,v_{\,i}^{\,2}\,(\,x\,)\,\left<\,\Lambda_{i}\,(\,x\,)\,P_{F_{\,i}\,(\,x\,)}\,U\,f,\, \Lambda_{i}\,(\,x\,)\,P_{F_{\,i}\,(\,x\,)}\,T\,f\,\right>\,d\mu_{x},\]
\[\Omega \,=\, \int\limits_{\,Y}\,v_{\,k}^{\,2}\,(\,x\,)\,\left<\,\Lambda_{i}\,(\,x\,)\,P_{F_{\,k}\,(\,x\,)}\,U\,f,\, \Lambda_{k}\,(\,x\,)\,P_{F_{\,k}\,(\,x\,)}\,T\,f\,\right>\,d\mu_{x}.\]
Then the family \,$\left\{\,\left(\,F_{\,i}\,(\,x\,),\, \Lambda_{i}\,(\,x\,),\, v_{\,i}\,(\,x\,)\,\right)\,\right\}_{x \,\in\, X,\, i \,\in\, [\,m\,]}$\, is W. C. C. K. G. F. F. for \,$H$\, with universal bounds \,$\dfrac{A}{(\,m \,-\, 1\,)\,(\,N \,+\, 1\,) \,+\, 1}$\, and \,$B$, where \,$A \,=\, \sum\limits_{i \,\in\, [\,m\,]}\,A_{i}$\, and \,$B \,=\, \sum\limits_{i \,\in\, [\,m\,]}\,B_{i}$.    
\end{theorem}

\begin{proof}
Let \,$\left\{\,\sigma_{i}\,\right\}_{i \,\in\, [\,m\,]}$\, be a partition of \,$X$.\,Then for \,$f \,\in\, H$, we have
\begin{align*}
&\sum\limits_{i \,\in\, [\,m\,]}\,A_{i}\,\left\|\,K^{\,\ast}\,f\,\right\|^{\,2} \,\leq\, \sum\limits_{i \,\in\, [\,m\,]}\,\int\limits_{\,X}\,v_{\,i}^{\,2}\,(\,x\,)\,\left<\,\Lambda_{i}\,(\,x\,)\,P_{F_{\,i}\,(\,x\,)}\,U\,f,\, \Lambda_{i}\,(\,x\,)\,P_{F_{\,i}\,(\,x\,)}\,T\,f\,\right>\,d\mu_{x}\\
&=\,\sum\limits_{i \,\in\, [\,m\,]}\,\sum\limits_{k \,\in\, [\,m\,]}\,\int\limits_{\,\sigma_{k}}\,v_{\,i}^{\,2}\,(\,x\,)\,\left<\,\Lambda_{i}\,(\,x\,)\,P_{F_{\,i}\,(\,x\,)}\,U\,f,\, \Lambda_{i}\,(\,x\,)\,P_{F_{\,i}\,(\,x\,)}\,T\,f\,\right>\,d\mu_{x} \\
&\leq\, \sum\limits_{i \,\in\, [\,m\,]}\,\bigg[\,\int\limits_{\,\sigma_{i}}\,v_{\,i}^{\,2}\,(\,x\,)\,\left<\,\Lambda_{i}\,(\,x\,)\,P_{F_{\,i}\,(\,x\,)}\,U\,f,\, \Lambda_{i}\,(\,x\,)\,P_{F_{\,i}\,(\,x\,)}\,T\,f\,\right>\,d\mu_{x} \,+\\ 
&\hspace{1cm}+\,\sum\limits_{k \,\in\, [\,m\,],\, k \,\neq\, i}\,\int\limits_{\,\sigma_{k}}\,\left<\,\Gamma_{i,\,k}\,U\,f,\, \Gamma_{i,\,k}\,T\,f\,\right>\,d\mu_{x} \,+\, \\
&\hspace{.5cm}\sum\limits_{k \,\in\, [\,m\,],\, k \,\neq\, i}\,\int\limits_{\,\sigma_{k}}\,v_{\,k}^{\,2}\,(\,x\,)\,\left<\,\Lambda_{k}\,(\,x\,)\,P_{F_{\,k}\,(\,x\,)}\,U\,f,\, \Lambda_{k}\,(\,x\,)\,P_{F_{\,k}\,(\,x\,)}\,T\,f\,\right>\,d\mu_{x}\,\bigg],\\
&\hspace{2cm}\Gamma_{i,\,k} \,=\, v_{\,i}^{\,2}\,(\,x\,)\,\Lambda_{i}\,(\,x\,)\,P_{F_{\,i}\,(\,x\,)}\,-\, v_{\,k}^{\,2}\,(\,x\,)\,\Lambda_{k}\,(\,x\,)\,P_{F_{\,k}\,(\,x\,)}\\
&\leq\, \sum\limits_{i \,\in\, [\,m\,]}\,\bigg[\,\int\limits_{\,\sigma_{i}}\,v_{\,i}^{\,2}\,(\,x\,)\,\left<\,\Lambda_{i}\,(\,x\,)\,P_{F_{\,i}\,(\,x\,)}\,U\,f,\, \Lambda_{i}\,(\,x\,)\,P_{F_{\,i}\,(\,x\,)}\,T\,f\,\right>\,d\mu_{x} \,+\\ 
&\hspace{.5cm}\sum\limits_{k \,\in\, [\,m\,],\, k \,\neq\, i}\,(\,N \,+\, 1\,)\,\int\limits_{\,\sigma_{k}}\,v_{\,k}^{\,2}\,(\,x\,)\,\left<\,\Lambda_{k}\,(\,x\,)\,P_{F_{\,k}\,(\,x\,)}\,U\,f,\, \Lambda_{k}\,(\,x\,)\,P_{F_{\,k}\,(\,x\,)}\,T\,f\,\right>\,d\mu_{x}\,\bigg],\\
&=\, \left\{\,(\,m \,-\, 1\,)\,(\,N \,+\, 1\,) \,+\, 1\,\right\}\,\sum\limits_{i \,\in\, [\,m\,]}\,\int\limits_{\,\sigma_{i}}\,v_{\,i}^{\,2}\,(\,x\,)\,\left<\,\Lambda_{i}\,(\,x\,)\,P_{F_{\,i}\,(\,x\,)}\,U\,f,\, \Lambda_{i}\,(\,x\,)\,P_{F_{\,i}\,(\,x\,)}\,T\,f\,\right>\,d\mu_{x}. 
\end{align*}
Thus, for each \,$f \,\in\, H$, we have
\begin{align*}
&\dfrac{A}{(\,m \,-\, 1\,)\,(\,N \,+\, 1\,) \,+\, 1}\,\left\|\,K^{\,\ast}\,f\,\right\|^{\,2}\\
&\leq\,\sum\limits_{i \,\in\, [\,m\,]}\,\int\limits_{\,\sigma_{i}}\,v_{\,i}^{\,2}\,(\,x\,)\,\left<\,\Lambda_{i}\,(\,x\,)\,P_{F_{\,i}\,(\,x\,)}\,U\,f,\, \Lambda_{i}\,(\,x\,)\,P_{F_{\,i}\,(\,x\,)}\,T\,f\,\right>\,d\mu_{x} \,\leq\, B\,\|\,f\,\|^{\,2}.
\end{align*}
This completes the proof.
\end{proof}

\section{Perturbation of woven continuous controlled $g$-fusion frame}

\smallskip\hspace{.6 cm}In frame theory, one of the most important problem is the stability of frame under some perturbation.\,P. Casazza and Chirstensen \cite{CC} have been generalized the Paley-Wiener perturbation theorem to perturbation of frame in Hilbert space.\,P. Ghosh and T. K. Samanta have studied perturbation of dual \,$g$-fusion frame and continuous controlled $g$-fusion frame in \cite{P, GG}.\,In this section, we will see that under some small perturbations, continuous controlled $K$-$g$-fusion frames constitute woven continuous controlled \,$K$-$g$-fusion frame.

\begin{theorem}
Let the families \,$\Lambda \,=\, \left\{\,\left(\,F\,(\,x\,),\, \Lambda\,(\,x\,),\, v\,(\,x\,)\,\right)\,\right\}_{x \,\in\, X}$\, and \,$\Gamma \,=\, \left\{\,\left(\,G\,(\,x\,),\, \Gamma\,(\,x\,),\, v\,(\,x\,)\,\right)\,\right\}_{x \,\in\, X}$\, be a continuous \,$(\,T,\,U\,)$-controlled \,$K$-$g$-fusion frame for \,$H$\, with bounds \,$A,\,B$\, and \,$C,\,D$, respectively.\,Suppose that there exist non-negative constants \,$\lambda_{1},\,\lambda_{2}$\, and \,$\mu$\, with \,$0 \,<\, \lambda_{\,1} \,<\, 1$, \,$\mu \,<\, \left(\,1 \,-\, \lambda_{\,1}\,\right)\,A \,-\, \lambda_{\,2}\,B$\, such that for each \,$f \,\in\, H$, we have 
\begin{align*}
&0\,\leq\\
&\int\limits_{\,X}\,v^{\,2}\,(\,x\,)\,\left<\,T^{\,\ast}\,\left(\,P_{\,F\,(\,x\,)}\,\Lambda\,(\,x\,)^{\,\ast}\,\Lambda\,(\,x\,)\,P_{\,F\,(\,x\,)}\,-\, P_{\,G\,(\,x\,)}\,\Gamma\,(\,x\,)^{\,\ast}\,\Gamma\,(\,x\,)\,P_{\,G\,(\,x\,)}\,\right)\,U\,f,\, f\,\right>\,d\mu_{x}\\
&\leq\,\lambda_{1}\,\int\limits_{\,X}\,v^{\,2}\,(\,x\,)\,\left<\,\Lambda\,(\,x\,)\,P_{\,F\,(\,x\,)}\,U\,f,\, \Lambda\,(\,x\,)\,P_{\,F\,(\,x\,)}\,T\,f\,\right>\,d\mu_{x} \,+\,\\
&\hspace{1cm}\lambda_{2}\,\int\limits_{\,X}\,v^{\,2}\,(\,x\,)\,\left<\,\Gamma\,(\,x\,)\,P_{\,G\,(\,x\,)}\,U\,f,\, \Gamma\,(\,x\,)\,P_{\,G\,(\,x\,)}\,T\,f\,\right>\,d\mu_{x} \,+\, \mu\,\left\|\,K^{\,\ast}\,f\,\right\|^{\,2}.   
\end{align*}
Then \,$\Lambda$\, and \,$\Gamma$\, are W. C. C. K. G. F. F. for \,$H$. 
\end{theorem}

\begin{proof}
Let \,$\sigma$\, be a partition of \,$X$.\,Now, for each \,$f \,\in\, H$, we have
\begin{align*}
& \int\limits_{\,\sigma}\,v^{\,2}\,(\,x\,)\,\left<\,\Lambda\,(\,x\,)\,P_{\,F\,(\,x\,)}\,U\,f,\, \Lambda\,(\,x\,)\,P_{\,F\,(\,x\,)}\,T\,f\,\right>\,d\mu_{x} \,+\,\\
&\hspace{1cm}\int\limits_{\,\sigma^{\,c}}\,v^{\,2}\,(\,x\,)\,\left<\,\Gamma\,(\,x\,)\,P_{\,G\,(\,x\,)}\,U\,f,\, \Gamma\,(\,x\,)\,P_{\,G\,(\,x\,)}\,T\,f\,\right>\,d\mu_{x}\\
&\geq\,\int\limits_{\,\sigma}\,v^{\,2}\,(\,x\,)\,\left<\,\Lambda\,(\,x\,)\,P_{\,F\,(\,x\,)}\,U\,f,\, \Lambda\,(\,x\,)\,P_{\,F\,(\,x\,)}\,T\,f\,\right>\,d\mu_{x}\,-\\
&\int\limits_{\sigma^{\,c}}\,v^{\,2}\,(\,x\,)\left<\,T^{\,\ast}\,\left(\,P_{\,F\,(\,x\,)}\,\Lambda\,(\,x\,)^{\,\ast}\,\Lambda\,(\,x\,)\,P_{\,F\,(\,x\,)}\,-\, P_{\,G\,(\,x\,)}\,\Gamma\,(\,x\,)^{\,\ast}\,\Gamma\,(\,x\,)\,P_{\,G\,(\,x\,)}\,\right)\,U\,f,\, f\,\right>\,d\mu_{x}\\
&\hspace{1cm}+\,\int\limits_{\,\sigma^{\,c}}\,v^{\,2}\,(\,x\,)\,\left<\,\Lambda\,(\,x\,)\,P_{\,F\,(\,x\,)}\,U\,f,\, \Lambda\,(\,x\,)\,P_{\,F\,(\,x\,)}\,T\,f\,\right>\,d\mu_{x}\\
&\geq\,\int\limits_{\,X}\,v^{\,2}\,(\,x\,)\,\left<\,\Lambda\,(\,x\,)\,P_{\,F\,(\,x\,)}\,U\,f,\, \Lambda\,(\,x\,)\,P_{\,F\,(\,x\,)}\,T\,f\,\right>\,d\mu_{x}\,-\\
&\int\limits_{\,X}\,v^{\,2}\,(\,x\,)\left<\,T^{\,\ast}\,\left(\,P_{\,F\,(\,x\,)}\,\Lambda\,(\,x\,)^{\,\ast}\,\Lambda\,(\,x\,)\,P_{\,F\,(\,x\,)}\,-\, P_{\,G\,(\,x\,)}\,\Gamma\,(\,x\,)^{\,\ast}\,\Gamma\,(\,x\,)\,P_{\,G\,(\,x\,)}\,\right)\,U\,f,\, f\,\right>\,d\mu_{x}\\
&\geq\,\left(\,1 \,-\, \lambda_{\,1}\,\right)\,\int\limits_{\,X}\,v^{\,2}\,(\,x\,)\,\left<\,\Lambda\,(\,x\,)\,P_{\,F\,(\,x\,)}\,U\,f,\, \Lambda\,(\,x\,)\,P_{\,F\,(\,x\,)}\,T\,f\,\right>\,d\mu_{x}\,-\\
&\hspace{1cm}\lambda_{2}\,\int\limits_{\,X}\,v^{\,2}\,(\,x\,)\,\left<\,\Gamma\,(\,x\,)\,P_{\,G\,(\,x\,)}\,U\,f,\, \Gamma\,(\,x\,)\,P_{\,G\,(\,x\,)}\,T\,f\,\right>\,d\mu_{x} \,-\, \mu\,\left\|\,K^{\,\ast}\,f\,\right\|^{\,2}\\
&\geq\,\bigg[\,\left(\,1 \,-\, \lambda_{\,1}\,\right)\,A \,-\, \lambda_{\,2}\,B \,-\, \mu\,\bigg]\,\left\|\,K^{\,\ast}\,f\,\right\|^{\,2}.
\end{align*}
On the other hand,
\begin{align*}
& \int\limits_{\,\sigma}\,v^{\,2}\,(\,x\,)\,\left<\,\Lambda\,(\,x\,)\,P_{\,F\,(\,x\,)}\,U\,f,\, \Lambda\,(\,x\,)\,P_{\,F\,(\,x\,)}\,T\,f\,\right>\,d\mu_{x} \,+\,\\
&\hspace{1cm}\int\limits_{\,\sigma^{\,c}}\,v^{\,2}\,(\,x\,)\,\left<\,\Gamma\,(\,x\,)\,P_{\,G\,(\,x\,)}\,U\,f,\, \Gamma\,(\,x\,)\,P_{\,G\,(\,x\,)}\,T\,f\,\right>\,d\mu_{x}\\
&\leq\, \int\limits_{\,X}\,v^{\,2}\,(\,x\,)\,\left<\,\Lambda\,(\,x\,)\,P_{\,F\,(\,x\,)}\,U\,f,\, \Lambda\,(\,x\,)\,P_{\,F\,(\,x\,)}\,T\,f\,\right>\,d\mu_{x} \,+\,\\
&\hspace{1cm}\int\limits_{\,X}\,v^{\,2}\,(\,x\,)\,\left<\,\Gamma\,(\,x\,)\,P_{\,G\,(\,x\,)}\,U\,f,\, \Gamma\,(\,x\,)\,P_{\,G\,(\,x\,)}\,T\,f\,\right>\,d\mu_{x}\\
&\leq\, (\,B \,+\, D\,)\,\|\,f\,\|^{\,2}.
\end{align*}
This completes the proof.
\end{proof}

\end{document}